\documentclass[10pt]{amsart}
\usepackage{latexsym, amsmath, amssymb}

\newcommand{\newsection}[1]
{\section{#1}\setcounter{theorem}{0} \setcounter{equation}{0}
\par\noindent}

\newtheorem{theorem}{Theorem}

\newtheorem{lemma}[theorem]{Lemma}

\newtheorem{proposition}[theorem]{Proposition}

\newcommand{\cd}{\, \cdot\, }

\newcommand{\R}{{\mathbb R}}
\newcommand{\ext}{{\R^3\backslash\mathcal{K}}}

\newcommand{\bdy}{{\partial\mathcal{K}}}

\newcommand{\ang}{{\not\negmedspace\nabla}}
\newcommand{\tr}{{\text{tr}}}

\begin{document}
\title[Null form wave equations in exterior domains]
{Global existence of null-form wave equations in exterior domains}
\thanks{The authors were supported in part by the NSF.  The second author was 
also supported by a fellowship from the Guggenheim Foundation}
\thanks{A portion of this work occurred while the authors were visiting the Mathematical
Sciences Research Institute, and the authors gratefully acknowledge the hospitality
and support of MSRI}

\author{Jason Metcalfe}
\address{Department of Mathematics, University of California, Berkeley, CA  94720-3840}
\email{metcalfe@math.berkeley.edu}

\author{Christopher D. Sogge}
\address{Department of Mathematics, Johns Hopkins University, Baltimore, MD  21218}
\email{sogge@jhu.edu}

\begin{abstract}
We provide a proof of global existence of solutions to quasilinear wave equations
satisfying the null condition in certain exterior domains.  In particular, our
proof does not require estimation of the fundamental solution for the free wave equation.  
We instead rely upon a class of Keel-Smith-Sogge estimates for the perturbed 
wave equation.  Using this, a notable simplification is made as compared to previous works
concerning wave equations in exterior domains: one no longer needs to distinguish
the scaling vector field from the other admissible vector fields.
\end{abstract}

\maketitle
\section{Introduction}

Inspired by the approach of Sideris and Tu \cite{Si3} in the boundaryless case and the application of
such techniques in Sideris \cite{Si}, we prove global existence for multi-speed systems
of quasilinear wave equations satisfying the null condition in certain exterior domains without
using estimation of the fundamental solution for the free wave equation.  To do so, we use 
a Keel-Smith-Sogge estimate for the perturbed equation established previously by the authors in
\cite{MS3}.  In the previous works on nonlinear wave equations in exterior domains, \cite{KSS3}, \cite{MS1,MS2},
and \cite{MNS1, MNS2}, it was necessary to use estimates that involved relatively few occurrences of the
scaling vector field $L=t\partial_t+r\partial_r$.  A notable innovation in the new approach allows us to
no longer distinguish between the scaling vector field $L$ and the other ``admissible'' invariant vector
fields $Z=\{\Omega_{ij}=x_i\partial_j - x_j\partial_i, \partial_k\,:\, 1\le i,j\le 3,\, 0\le k\le 3\}$.  This is
accomplished by introducing modified vector fields that preserve the boundary condition, as was in part done in \cite{MS1}.

The main existence result is an analog of the classic results of Christodoulou \cite{Christodoulou} and Klainerman 
\cite{K} in the
boundaryless case.  In the multiple speed boundaryless case, related results were established, e.g., by
Sideris and Tu \cite{Si3}, Sogge \cite{So2}, Agemi and Yokoyama \cite{AY}, 
Kubota and Yokoyama \cite{KY}, Hidano \cite{Hid}, Yokoyama \cite{Y}, and  Katayama \cite{Kata, Kata2, Kata3}.  
In exterior domains, null form quasilinear wave equations
were previously studied by Keel, Smith, and Sogge \cite{KSS}, the authors \cite{MS1}, and Metcalfe, Nakamura,
and Sogge \cite{MNS1, MNS2}.  We note that the main theorem of this paper was previously established
(in a more general context) in \cite{MNS1}.  We, however, believe that the new techniques are of independent interest,
and we are hopeful about their potential use in applications.

Our proof uses Klainerman's method of commuting vector fields \cite{K} as was adapted to exterior domains
by Keel, Smith, and Sogge \cite{KSS3}.  In particular, we restrict to the class of ``admissible'' vector
fields that was mentioned above.  Notably absent in this set are the hyperbolic rotations
$\Omega_{0j}=t\partial_j + x_j \partial_t$ which do not seem appropriate for problems in exterior domains as they
have unbounded normal component on the boundary.  Moreover, in the multiple speed setting, these vector fields
have an associated speed and only commute with the d'Alembertian of the same speed.

This approach relies upon a weighted space-time $L^2$ estimate, which will be referred to as the Keel-Smith-Sogge
estimate or KSS estimate.  Such estimates were established in \cite{KSS2} where they were first used to
show long-time existence of solutions to nonlinear equations.  With this estimate, existence is established using
$O(1/|x|)$ decay rather than the more standard, but quite difficult to prove when there is boundary, $O(1/t)$ decay.
An earlier, related estimate is due to Strauss
\cite{Strauss} (Lemma 3).  The proof in \cite{KSS2} is easily modified to establish these bounds in all dimensions
$n\ge 3$ as is done in Metcalfe \cite{Met} and Hidano and Yokoyama 
\cite{HidY} and has been used in, e.g., \cite{HidY}, \cite{Met},
and \cite{MS2} to study nonlinear equations.  Recently, using techniques of Rodnianski \cite{Sterb},
the authors \cite{MS3} have established an analogous estimate for the perturbed equation.  This inequality 
is essential to the approach presented in this article.  It is worth noting that Alinhac \cite{Alinhac2} 
simultaneously obtained
a related KSS-type estimate for the perturbed equation and for wave equations on curved backgrounds.  The assumptions
made on the perturbation in \cite{Alinhac2} are, however, not as favorable in the current setting.

The KSS estimate and energy estimates will be coupled with some well-known decay estimates in order to get global existence.  
These decay estimates are variants of those of Klainerman and Sideris \cite{KS} and are known to be rather widely applicable,
including e.g. to the equations of elasticity.  This will, of course, be used in combination with the extra
decay afforded to us by the null condition.  

The major innovations of this paper regard the variable-coefficient KSS estimate.  First, we expand upon the proof
in \cite{MS3} and carefully prove the KSS estimates for the perturbed wave equation in the multiple speed setting.
Having such estimates for the perturbed equation allows us to apply the KSS estimates even for terms of the highest order.
This was not previously possible for quasilinear equations as there was a loss of regularity resulting from the occurrence
of second derivatives.  As such, we may now prove global existence using only energy methods.  In particular, the
decay estimates that we shall employ do not require direct estimation of the fundamental solution of the linear
wave equation, and in particular, such estimates are known to hold for some related applications.  We note, e.g., that 
the only obstacle to using similar techniques to study the equations of isotropic elasticity in exterior domains is
deriving the existence of a KSS-type estimate for the perturbed linearized equations.  This, however,
is more delicate than in the current setting due to the off-diagonal terms and is currently an open problem of interest.

A key obstacle to using only energy methods in previous studies was the boundary terms that arise in 
the Klainerman-Sideris estimates.  In particular, there is a term localized near the boundary which has
significantly less decay in $t$.  To alleviate this, we shall use the additional decay in $x$, which was 
largely ignored in previous works, from the Klainerman-Sideris estimates.  When combined with the KSS estimates,
this decay permits us the necessary control over the boundary term.

As the coefficients of the scaling vector field can be large in an arbitrarily small neighborhood of the obstacle,
previous works in exterior domains using the adapted method of commuting vector fields use estimates that
required relatively few occurrences of $L$, and during the proofs of long-time existences, the scaling vector fields
must be carefully tracked.  This, at best, complicated these arguments. 
Using the variable coefficient KSS estimates, we show that it
is not strictly necessary to differentiate $L$ from the other admissible vector fields. 
In particular, we use a modified scaling vector field which preserves
the Dirichlet boundary conditions as in \cite{MS1}, but instead of proving boundary term estimates, such as
\cite[Lemma 2.9]{MS1}, we are now able to control the resulting commutators using the variable coefficient KSS
estimates.

Finally, we mention that our proof, unlike many of the previous works in exterior domains,
 does not directly use the decay of local energy, such as that of Lax,
Morawetz, and Phillips \cite{LMP}.  Our hypotheses on the 
obstacle are, however, sufficient to guarantee said decay.  It is conceivable that the techniques
contained herein
could be important in other applications where the rate of decay of local energy is slower, say e.g. in even dimensions.

We now more precisely describe the main result of this article.  We begin by fixing a bounded obstacle
$\mathcal{K}\subset\R^3$ with smooth boundary.  Moreover, we shall assume that $\mathcal{K}$ is star-shaped with respect
to the origin.  As we shall see, scaling will allow us to assume without loss that $\mathcal{K}\subset \{|x|\le 1\}$, and
this assumption is made throughout.  
The star-shapedness assumption is used to see 
that certain boundary terms in our energy estimates
and KSS estimates have a favorable sign.  This is reminiscent of arguments from Morawetz \cite{M}.

In the exterior of $\mathcal{K}$, we shall study systems of quasilinear wave equations of the form
\begin{equation}
\label{i.1}
\begin{cases}
\partial_t^2 u^I - c_I^2 \Delta u^I = B^{IJ,\alpha\beta}_{K,\gamma} \partial_\gamma u^K \partial_\alpha\partial_\beta u^J,
\quad (t,x)\in \R_+\times\ext,\quad I=1,2,\dots,D,\\
u(t,\cd)|_{\bdy} =0,\\
u(0,\cd)=f,\quad \partial_tu(0,\cd)=g.
\end{cases}
\end{equation}
Here and throughout we use the Einstein summation convention.  Repeated Greek indices $\alpha,\beta,\gamma$, 
and $\delta$ are implicitly summed from $0$ to $3$.  
  Repeated lowercase Latin indices $a,b$ are summed from $1$
to $3$, and repeated uppercase Latin indices $I,J,K$ are summed from $1$ to $D$. 
In the sequel, we will use $\Box=(\Box_{c_1},\dots,\Box_{c_D})$ to denote the vector-valued d'Alembertian,
where $\Box_{c_I}=\partial_t^2 - c_I^2\Delta$.  For simplicity, we shall study the nonrelativistic case where
$$c_1>\dots>c_D>0.$$
Straightforward modifications will allow for repeated wave speeds.  

The $B_{K,\gamma}^{IJ,\alpha\beta}$ in \eqref{i.1} are real constants satisfying the symmetry conditions
\begin{equation}\label{i.2}
B^{IJ,\alpha\beta}_{K,\gamma}=B^{JI,\alpha\beta}_{K,\gamma}=B^{IJ,\beta\alpha}_{K\gamma}.
\end{equation}
In order to get global existence, we shall also need to assume that \eqref{i.1} satisfies the null condition.
In the nonrelativistic case, this says that the self-interactions among the quasilinear terms 
satisfy the standard null condition.  That is,
\begin{equation}
\label{i.3}
B^{JJ,\alpha\beta}_{J,\gamma}\xi_\alpha\xi_\beta\xi_\gamma = 0,\quad\text{whenever } \frac{\xi_0^2}{c_J^2}-
\xi_1^2-\xi_2^2-\xi_3^2 =0, \quad J=1,\dots,D.
\end{equation}

To solve \eqref{i.1}, the data must be assumed to satisfy the relevant compatibility conditions.  Letting
$J_ku=\{\partial_x^\alpha u\, :\, 0\le |\alpha|\le k\}$, we know that for a fixed $m$ and formal $H^m$ solution
$u$ of \eqref{i.1}, we can write $\partial_t^k u(0,\cd)=\psi_k (J_kf, J_{k-1}g)$, $0\le k\le m$, for 
compatibility functions $\psi_k$ depending on the nonlinearity, $J_kf$, and $J_{k-1}g$.  For $(f,g)\in H^m\times
H^{m-1}$, the compatibility condition simply requires that the $\psi_k$ vanish on $\bdy$ for $0\le k\le m-1$.
For smooth $(f,g)$, we say that the compatibility condition is satisfied to infinite order if this vanishing
condition holds for all $m$.  We refer the reader to \cite{KSS} for a more thorough exposition on these
compatibility conditions.

We are now prepared to state our main theorem.
\begin{theorem}
\label{theoremi.1}
Let $\mathcal{K}$ be a fixed compact obstacle with smooth boundary that is star-shaped with 
respect to the origin.  Assume that the $B_{K,\gamma}^{IJ,\alpha\beta}$ are as above.  Then, there 
is a constant $\varepsilon_0>0$ and an integer constant $N>0$ so that if the data $(f,g)\in
C^\infty(\ext)$ satisfy the compatibility condition to infinite order and the smallness
condition
\begin{equation}
\label{i.4}
\sum_{|\alpha|\le N} \|\langle x\rangle^{|\alpha|}\partial^\alpha_x \partial_x f\|_2
+\sum_{|\alpha|\le N} \|\langle x\rangle^{|\alpha|} \partial^\alpha_x g\|_2\le \varepsilon_0,
\end{equation}
then \eqref{i.1} has a unique global solution $u\in C^\infty([0,\infty)\times \ext)$.
\end{theorem}

This paper is organized as follows.  In the next section, we provide derivations of the energy estimates
and KSS estimates for the perturbed wave equations.  We also show that an appropriate variant of these
holds when $u$ is replaced by $\Gamma^\mu u$.  Here $\Gamma=\{L,Z\}$ is the set of admissible vector
fields.  As mentioned previously, the fact that we no longer need to distinguish between $L$ and $Z$ is a
significant innovation in this paper.  In the third section, we present the decay estimates that we will require.
These are fairly well-known, but in the interest of making this paper somewhat self-contained, the proofs
are sketched.  In the last section, we prove the main result, Theorem \ref{theoremi.1}.

\bigskip
\newsection{Energy estimates and Keel-Smith-Sogge estimates}

%

In this section, we establish the energy and KSS estimates for the perturbed wave equation 
that we shall require in the sequel.  We must take care to insure that our estimates will not destroy the 
null structure.

We will be concerned with solutions $u^I\in C^\infty(\R_+\times\ext)$ of the Dirichlet-wave equation
\begin{equation}
\label{e.1}
\begin{cases}
(\Box_h u)^I = F^I,\\
u|_{\bdy}=0
\end{cases}
\end{equation}
where
\begin{equation}
\label{e.2}
(\Box_h u)^I = (\partial^2_t - c_I^2 \Delta)u^I + 
\sum_{J=1}^D \sum_{\alpha,\beta=0}^3 h^{IJ,\alpha\beta}(t,x)\partial_\alpha\partial_\beta u^J.
\end{equation}
We shall assume that the $h^{IJ,\alpha\beta}$ satisfy the symmetry conditions
\begin{equation}\label{e.3}
h^{IJ,\alpha\beta}=h^{JI,\alpha\beta}=h^{IJ,\beta\alpha},
\end{equation}
as well as the size condition
\begin{equation}
\label{e.4}
|h|=\sum_{I,J=1}^D \sum_{\alpha,\beta=0}^3 |h^{\alpha\beta}(t,x)|\le \delta \ll 1.
\end{equation}

We denote $u=(u^1, \dots, u^D)$.  
Here, we are working in the Euclidean metric, and indices are raised with this metric.  

We will need to define the full energy-momentum tensor associated to \eqref{e.1}.  To begin,
let
\begin{multline}
\label{e.5}
Q_{0 \beta}[u] = \partial_0 u^I \partial_\beta u^I -\frac{1}{2}\delta_{0\beta}\Bigl[
|\partial_0 u|^2 - c_I^2 |\nabla_x u^I|^2\Bigr] \\+ \delta_{0\gamma} h^{IJ,\gamma\delta} \partial_\delta u^J
\partial_\beta u^I - \frac{1}{2}\delta_{0\beta} h^{IJ,\gamma\delta} \partial_\gamma u^J \partial_\delta u^I,
\end{multline}
and
\begin{multline}
\label{e.6}
Q_{\alpha\beta}[u]=-c_I^2\partial_\alpha u^I \partial_\beta u^I - \frac{1}{2}\delta_{\alpha\beta}\Bigl[|\partial_0 u|^2
-c_I^2 |\nabla_x u^I|^2\Bigr] \\+ \delta_{\alpha\gamma} h^{IJ,\gamma\delta}\partial_\delta u^J \partial_\beta u^I
-\frac{1}{2}\delta_{\alpha\beta} h^{IJ,\gamma\delta}\partial_\gamma u^J \partial_\delta u^I,\quad \alpha=1,2,3.
\end{multline}
An elementary calculation yields
\begin{equation}
\label{e.7}
D^\alpha Q_{\alpha\beta}[u]= \partial_\beta u^I (\Box_h u)^I 
+ (\partial_\gamma h^{IJ,\gamma\delta}) \partial_\delta u^J \partial_\beta u^I -\frac{1}{2}(\partial_\beta h^{IJ,\gamma\delta})
\partial_\gamma u^J \partial_\delta u^I.
\end{equation}

\subsection{Energy estimate}

From \eqref{e.7}, we are quickly able to obtain the well-known energy estimate for the perturbed wave
equation.
\begin{proposition}
\label{prope.1}
Assume that $\mathcal{K}$ is a bounded obstacle with $C^1$-boundary.
Assume also that the perturbation terms are as above.  Suppose that $u\in C^\infty$ solves
\eqref{e.1} and for every $t$, $u(t,x)=0$ for large $x$.  Then, 
\begin{multline}
\label{e.8} \|u'(t,\cd)\|^2_2 \lesssim \|u'(0,\cd)\|^2_2 + \int_0^t \int_\ext |(\Box_h u)^I \partial_0 u^I|
\:dx\:ds \\+ \int_0^t \int_\ext \Bigl[|(\partial_\gamma h^{IJ,\gamma\delta}) \partial_\delta u^J 
\partial_0 u^I| + |(\partial_0 h^{IJ,\gamma\delta})
\partial_\gamma u^J \partial_\delta u^I|\Bigr]\:dx\:ds.
\end{multline}
\end{proposition}

Here $u'=(\partial_t u,\nabla_x u)$ is used to denote the full space-time gradient.

Indeed, we need only examine the $\beta=0$ components of \eqref{e.5} and \eqref{e.6}.  Integrating
\eqref{e.7} over $S_t=[0,t]\times \ext$, it immediately follows that
\begin{multline}\label{e.9}
\int_\ext Q_{00}[u](t,\cd)\:dx = \int_\ext Q_{00}[u](0,\cd)\:dx + \int_0^t \int_\ext
\partial_0 u^I (\Box_h u)^I \:dx\:ds \\+ \int_0^t \int_\ext (\partial_\gamma h^{IJ,\gamma\delta})\partial_\delta
u^J \partial_0 u^I \:dx\:ds -\frac{1}{2} \int_0^t \int_\ext (\partial_0 h^{IJ,\gamma\delta}) \partial_\gamma u^J 
\partial_\delta u^I.
\end{multline}
Here, we have used the fact that $\partial_t$ preserves the Dirichlet boundary condition.  Thus, the integrand 
of the boundary term that results in the application of the divergence theorem vanishes identically.
If $\delta$ in \eqref{e.4} is sufficiently small, it follows that
\begin{equation}
\label{e.10}
(2\max_I \{c_I^2,c_I^{-2}\})^{-1}|u'(t,x)|^2 \le Q_{00}[u](t,x)\le 2\max_I\{c_I^2,c_I^{-2}\} |u'(t,x)|^2.
\end{equation}
And, thus, \eqref{e.9} immediately yields \eqref{e.8}.

\subsection{Keel-Smith-Sogge estimates}
As mentioned previously, a key estimate that allows long time existence to be deduced from decay
in the spatial variables is a weighted mixed norm estimate of Keel, Smith, and Sogge \cite{KSS2}.  
In a different context, Rodnianski \cite{Sterb} proved a variant of the KSS estimate using energy
methods, and as was shown in \cite{MS3}, these methods are stable under small perturbations.  
The following proposition is essentially from \cite{MS3}.  Here, additional care is required
to preserve the null structure of the equation, and we also explicitly examine the multiple
speed system.

\begin{proposition}
\label{prope.2}
Suppose that $\mathcal{K}$ is a $C^1$, bounded, star-shaped obstacle as above.  Suppose,
further, that the perturbation terms $h^{IJ,\alpha\beta}$ are as above.  Then, if $u\in C^\infty$ 
solves \eqref{e.1} and for every $t$, $u(t,x)=0$ for large $x$, 
\begin{multline}
\label{e.11}
\|\langle x\rangle^{-1/2-} u'\|^2_{L^2_tL^2_x(S_t)} + (\log(2+t))^{-1} \|\langle x\rangle^{-1/2}
u'\|^2_{L^2_tL^2_x(S_t)}
\\\lesssim \|u'(0,\cd)\|^2_2 + \int_0^t \int_\ext \Bigl(|\partial_{t,x} u^I| + 
\frac{|u^I|}{r}\Bigr) |(\Box_h u)^I|\:dx\:ds
\\+\int_0^t \int_\ext \Bigl|(\partial_\gamma h^{IJ,\gamma\delta})\partial_\delta u^J\Bigr| \Bigl(|\partial_{t,r} u^I|
+ \frac{|u^I|}{r}\Bigr)\:dx\:ds
\\+ \int_0^t \int_\ext \Bigl|(\partial_{t,r}
h^{IJ,\gamma\delta})\partial_\gamma u^J \partial_\delta u^I\Bigr| 
\:dx\:ds
+\int_0^t \int_\ext \frac{|h|}{\langle x\rangle} |\nabla u|\Bigl(|\nabla u|+\frac{|u|}{r}\Bigr)\:dx\:ds
\end{multline}
for any $t\ge 0$.
\end{proposition}

In the proposition, we are using the notation $\langle x\rangle = \langle r\rangle = \sqrt{1+|x|^2}$.  We
are also using the notation $\langle x\rangle^{-1/2-}$ to indicate that the estimate holds with this 
weight replaced by $\langle x\rangle^{-1/2-\delta}$ for any $\delta>0$.  The implicit constant depends
on this $\delta$, and in practice, we will only require the estimate for a fixed, positive $\delta$.

Here, we contract the energy momentum tensor, \eqref{e.5} and \eqref{e.6}, with a radial vector
field $X=f(r)\partial_r$ which allows us to define the momentum density 
$$P_\alpha[u,X]=Q_{\alpha\beta}[u]X^\beta.$$
Computing the divergence, we have
\begin{multline*}
-D^\alpha P_\alpha[u,X]=-(\partial_r u^I)(\Box_h u)^I f(r)
-(\partial_\gamma h^{IJ,\gamma \delta}) \partial_\delta u^J \partial_r u^I f(r)
\\+\frac{1}{2} (\partial_r h^{IJ,\gamma\delta})\partial_\gamma u^J \partial_\delta u^I f(r)
\\+c_I^2 f'(r) (\partial_r u^I)^2 + c_I^2 \frac{f(r)}{r} |\ang u^I|^2 - \frac{1}{2}\tr\pi
\Bigl[-|\partial_0 u|^2 +c_I^2 |\nabla_x u^I|^2\Bigr]
\\-\frac{x_a}{r} h^{IJ,a\delta} \partial_\delta u^J \partial_r u^I f'(r) 
+\frac{x_a}{r}h^{IJ,a\delta}\partial_\delta u^J \partial_r u^I \frac{f(r)}{r} 
\\-h^{IJ,a\delta}\partial_a u^I \partial_\delta u^J \frac{f(r)}{r} 
+ \frac{1}{2}(\tr\pi)
h^{IJ,\gamma\delta} \partial_\gamma u^J \partial_\delta u^I.
\end{multline*}
Here $\pi$ denotes the deformation tensor of $X$ and, as can be checked,
\begin{equation}
\label{e.12}
\tr\pi = f'(r) + 2 \frac{f(r)}{r}.
\end{equation}

At this point, we define the modified momentum density
\begin{equation}
\label{e.13}
\tilde{P}_0[u,X]=P_0[u,X] + \frac{f(r)}{r} u^I \partial_0 u^I + \frac{f(r)}{r} h^{IJ,0\beta} u^I \partial_\beta u^J,
\end{equation}
\begin{multline}
\label{e.14}
\tilde{P}_\alpha[u,X]=P_\alpha[u,X]- c_I^2 \frac{f(r)}{r} u^I \partial_\alpha u^I + \frac{c_I^2}{2}
\partial_\alpha \Bigl(\frac{f(r)}{r}\Bigr)(u^I)^2 + \frac{f(r)}{r} h^{IJ,\alpha\beta} u^I \partial_\beta u^J,\quad
\\\alpha=1,2,3
\end{multline}
A tedious but elementary calculation yields
\begin{multline}
\label{e.15}
-D^\alpha \tilde{P}_\alpha[u,X]=
-\Bigl(\partial_r u^I +\frac{u^I}{r}\Bigr)(\Box_h u)^I f(r)
-(\partial_\gamma h^{IJ,\gamma \delta}) \partial_\delta u^J \Bigl(\partial_r u^I+\frac{u^I}{r}\Bigr) f(r)
\\+\frac{1}{2} (\partial_r h^{IJ,\gamma\delta})\partial_\gamma u^J \partial_\delta u^I f(r)
\\+c_I^2 f'(r) (\partial_r u^I)^2 + c_I^2 \frac{f(r)}{r} |\ang u^I|^2 - \frac{1}{2}f'(r)
\Bigl[-|\partial_0 u|^2 +c_I^2 |\nabla_x u^I|^2\Bigr]
\\-\frac{x_a}{r} h^{IJ,a\delta} \partial_\delta u^J \Bigl(\partial_r u^I+\frac{u^I}{r}\Bigr) f'(r) 
+\frac{x_a}{r}h^{IJ,a\delta}\partial_\delta u^J \Bigl(\partial_r u^I+\frac{u^I}{r}\Bigr) \frac{f(r)}{r} 
\\
-h^{IJ,a\delta}\partial_a u^I \partial_\delta u^J \frac{f(r)}{r} 
+\frac{1}{2}f'(r)
h^{IJ,\gamma\delta} \partial_\gamma u^J \partial_\delta u^I - \frac{c_I^2}{2}\Delta\Bigl(\frac{f(r)}{r}\Bigr)(u^I)^2
\end{multline}

Integrating both sides of \eqref{e.14} in a time strip $S_t$ gives
\begin{multline}
\label{e.16}
\int_\ext \tilde{P}_0[u,X](0)\:dx - \int_\ext \tilde{P}_0[u,X](t)\:dx +\int_0^t \int_{\bdy}
\tilde{P}_a[u,X]n^a\:d\sigma\:ds
\\=-\int_0^t\int_\ext D^\alpha \tilde{P}_\alpha[u,X]\:dx\:ds.
\end{multline}
Here $\vec{n}=(n^1,n^2,n^3)$ is the outward unit normal to $\mathcal{K}$, and $d\sigma$ is the surface
measure on $\bdy$.

At this point, as in \cite{Sterb} and \cite{MS3}, we choose
$$f(r)=\frac{r}{r+\rho}$$
for a positive constant $\rho$.
Notice, in particular, that we have $|f(r)|\lesssim 1$ and $|f'(r)|\lesssim \frac{1}{r}$.  Thus, it
follows that
\begin{equation}\label{e.17}
\begin{split}
\Bigl|\int_\ext \tilde{P}_0[u,X](0)\:dx\Bigr|&=\Bigl|\int_\ext \Bigl(\partial_t u(0,x)\partial_r u(0,x) f(r)
\\-f(r)&h^{IJ,0\delta}(0,x)\partial_\delta u^J(0,x) \partial_r u^I(0,x)
 + \frac{f(r)}{r}u^I(0,x) \partial_t u^I(0,x)
\\&\qquad\qquad\qquad\qquad\qquad+\frac{f(r)}{r}
h^{IJ,0\beta}(0,x) u^I(0,x)\partial_\beta u^J(0,x)\Bigr)\:dx\Bigr|
\\
&\lesssim \|u'(0,\cd)\|^2_2.
\end{split}\end{equation}
For the last inequality, we are applying the Schwarz inequality and a Hardy inequality.  We are also using
\eqref{e.4}.

A similar bound holds for $\tilde{P}_0[u,X](t)$, and thus, using the energy inequality \eqref{e.8},
\begin{equation}
\label{e.18}
\begin{split}
\Bigl|\int_\ext \tilde{P}_0[u,X](t)\:dx\Bigr| &\lesssim \|u'(t,\cd)\|^2_2\\
&\lesssim \|u'(0,\cd)\|^2_2 + \int_0^t\int_\ext |(\Box_h u)^I\partial_0 u^I|\:dx\:ds\\
& + \int_0^t \int_\ext \Bigl[|(\partial_\gamma h^{IJ,\gamma\delta})\partial_\delta u^J \partial_0 u^I|
+|(\partial_0 h^{IJ,\gamma\delta})\partial_\gamma u^J \partial_\delta u^I|\Bigr]\:dx\:ds.
\end{split}
\end{equation}

Since the Dirichlet boundary conditions allow us to write $\partial_a u^I = \partial_{\vec{n}} u^I n_a$ on 
$\bdy$ and since $\langle x, \vec{n}\rangle>0$ on $\bdy$ for star-shaped $\mathcal{K}$, for the spatial
boundary terms, we have
\begin{equation}
\label{e.19}
\int_0^t \int_\bdy \tilde{P}_a[u,X]n^a\:d\sigma\:ds\le -\frac{c_I^2}{4}\int_0^t \int_\bdy \frac{f(r)}{r}
(\partial_{\vec{n}} u^I)^2 \langle x,\vec{n}\rangle\:d\sigma\:ds\le 0.
\end{equation}
Here, we have also used the smallness of the perturbation, \eqref{e.4}.

If we use that $\Delta(f(r)/r)\le 0$ and \eqref{e.17}-\eqref{e.19} in \eqref{e.16}, we see
that
\begin{multline*}
\int_0^t \int_\ext f'(r) c_I^2 (\partial_r u^I)^2 + \frac{f(r)}{r} c_I^2 |\ang u^I|^2 -\frac{1}{2}
f'(r)\Bigl[-|\partial_0 u|^2 + c_I^2 |\nabla_x u^I|^2\Bigr]\:dx\:ds
\\\lesssim \|u'(0,\cd)\|^2_2 + \int_0^t \int_\ext \Bigl(|\partial_{t,x} u^I| + 
\frac{|u^I|}{r}\Bigr) |(\Box_h u)^I|\:dx\:ds
\\+\int_0^t \int_\ext \Bigl|(\partial_\gamma h^{IJ,\gamma\delta})\partial_\delta u^J\Bigr| \Bigl(|\partial_{t,r} u^I|
+ \frac{|u^I|}{r}\Bigr)\:dx\:ds
\\+ \int_0^t \int_\ext \Bigl|(\partial_{t,r}
h^{IJ,\gamma\delta})\partial_\gamma u^J \partial_\delta u^I\Bigr| 
\:dx\:ds
+\int_0^t \int_\ext \frac{|h|}{\langle x\rangle} |\nabla u|\Bigl(|\nabla u|+\frac{|u|}{r}\Bigr)\:dx\:ds.
\end{multline*}
Since $f'(r)\le\frac{f(r)}{r}$, this implies
\begin{multline}
\label{e.20}
\int_0^t \int_\ext f'(r) |\partial_t u|^2 + f'(r) c_I^2 (\partial_r u^I)^2 + \frac{f(r)}{r} c_I^2 |\ang u^I|^2 \:dx\:ds
\\\lesssim \|u'(0,\cd)\|^2_2 + \int_0^t \int_\ext \Bigl(|\partial_{t,x} u^I| + 
\frac{|u^I|}{r}\Bigr) |(\Box_h u)^I|\:dx\:ds
\\+\int_0^t \int_\ext \Bigl|(\partial_\gamma h^{IJ,\gamma\delta})\partial_\delta u^J\Bigr| \Bigl(|\partial_{t,r} u^I|
+ \frac{|u^I|}{r}\Bigr)\:dx\:ds
\\+ \int_0^t \int_\ext \Bigl|(\partial_{t,r}
h^{IJ,\gamma\delta})\partial_\gamma u^J \partial_\delta u^I\Bigr| 
\:dx\:ds
+\int_0^t \int_\ext \frac{|h|}{\langle x\rangle} |\nabla u|\Bigl(|\nabla u|+\frac{|u|}{r}\Bigr)\:dx\:ds.
\end{multline}

By choosing $\rho=1$ and $\rho=2^k$ for an integer $k\ge 1$ respectively, we see that
$$\int_0^t \int_{|x|\le 1} |u'|^2\:dx\:ds$$
and
$$\int_0^t \int_{2^{k-1}\le |x|\le 2^k} \frac{|u'|^2}{r}\:dx\:ds$$
are bounded by the right side of \eqref{e.20}.  If we sum these resulting estimates over $k\ge 1$, 
we see immediately that the bound for the first term in \eqref{e.11} holds.  The same argument 
yields the bound for the second term in the left of \eqref{e.11}.  Indeed, the estimate follows trivially
from \eqref{e.8} when the spatial norm is over $|x|\ge t$.  Thus, we need only sum over the $O(\log(2+t))$ 
choices of $k$ with $2^{k-1}\lesssim t$.

\subsection{Main $L^2$ estimate}

In this section, we show that higher order energy estimates also hold.  In particular, we show that
versions of \eqref{e.8} and \eqref{e.11} hold when $u$ is replaced by $\Gamma^\mu u$.  In order to do so,
we introduce modified vector fields that preserve the boundary condition.  This extends an idea
initiated in \cite{MS1}.

Notice that, by combining the main results (\eqref{e.8} and \eqref{e.11}) of the preceding sections, we have
\begin{multline}
\label{e.21}
\|\langle x\rangle^{-1/2-} u'\|^2_{L^2_tL^2_x(S_t)} + (\log(2+t))^{-1} \|\langle x\rangle^{-1/2}
u'\|^2_{L^2_tL^2_x(S_t)} + \|u'(t,\cd)\|^2_2
\\\lesssim \|u'(0,\cd)\|^2_2 + \int_0^t \int_\ext \Bigl(|\partial_{t,x} u^I| + 
\frac{|u^I|}{r}\Bigr) |F^I|\:dx\:ds
\\+\int_0^t \int_\ext \Bigl|(\partial_\gamma h^{IJ,\gamma\delta})\partial_\delta u^J\Bigr| \Bigl(|\partial_{t,r} u^I|
+ \frac{|u^I|}{r}\Bigr)\:dx\:ds
\\+ \int_0^t \int_\ext \Bigl|(\partial_{t,r}
h^{IJ,\gamma\delta})\partial_\gamma u^J \partial_\delta u^I\Bigr| 
\:dx\:ds
+\int_0^t \int_\ext \frac{|h|}{\langle x\rangle} |\nabla u|\Bigl(|\nabla u|+\frac{|u|}{r}\Bigr)\:dx\:ds
\end{multline}
when $u$ solves \eqref{e.1}.  Notice, in particular, that if $F$ vanishes for $|x|>2$, then
we can bound the second term in the right side by
$$\|u'\|_{L^2_tL^2_x(S_t\cap \{|x|\le 2\})}\|F\|_{L^2_tL^2_x(S_t)} \le \varepsilon \|\langle x
\rangle^{-1/2-} u'\|^2_{L^2_tL^2_x(S_t)} + C \|F\|^2_{L^2_tL^2_x(S_t)},$$
and, in this case, the first term on the right of this inequality can be bootstrapped.  
Here, we have used the fact that the Dirichlet
boundary condition allows us to control $u$ locally by $u'$.  We have also used that $0\in \mathcal{K}$,
and hence, $1/r$ is bounded on $\ext$.

Thus, it immediately follows that if $u$ is a solution to
\begin{equation}
\label{e.22}
\begin{cases}
\Box_h u = F+G\\
u|_\bdy=0,
\end{cases}
\end{equation}
and $G$ vanishes unless $|x|\le 2$, then
\begin{multline}
\label{e.23}
\|\langle x\rangle^{-1/2-} u'\|^2_{L^2_tL^2_x(S_t)} + (\log(2+t))^{-1} \|\langle x\rangle^{-1/2}
u'\|^2_{L^2_tL^2_x(S_t)} + \|u'(t,\cd)\|^2_2
\\\lesssim \|u'(0,\cd)\|^2_2 + \int_0^t \int_\ext \Bigl(|\partial_{t,x} u^I| + 
\frac{|u^I|}{r}\Bigr) |F^I|\:dx\:ds + \|G\|^2_{L^2_tL^2_x(S_t)}
\\+\int_0^t \int_\ext \Bigl|(\partial_\gamma h^{IJ,\gamma\delta})\partial_\delta u^J\Bigr| \Bigl(|\partial_{t,r} u^I|
+ \frac{|u^I|}{r}\Bigr)\:dx\:ds
\\+ \int_0^t \int_\ext \Bigl|(\partial_{t,r}
h^{IJ,\gamma\delta})\partial_\gamma u^J \partial_\delta u^I\Bigr| 
\:dx\:ds
+\int_0^t \int_\ext \frac{|h|}{\langle x\rangle} |\nabla u|\Bigl(|\nabla u|+\frac{|u|}{r}\Bigr)\:dx\:ds.
\end{multline}
We will use this as a base case for an induction argument to construct higher order energy estimates.

Since $\partial_t^j$ preserves the Dirichlet boundary condition, the estimate \eqref{e.23} holds with
$u$ replaced by $\partial_t^j u$.  Moreover, if we apply elliptic regularity (see, e.g., \cite{MS1} Lemma 2.3), 
it follows that
\begin{multline}
\label{e.24}
\sum_{|\mu|\le N} \|\partial^\mu u'\|^2_{L^2_tL^2_x(S_t\cap\{|x|\le 2\})} \lesssim \sum_{j\le N}
\|\partial_t^j u'(0,\cd)\|^2_2 \\+ \sum_{j,k\le N} \int_0^t \int_\ext \Bigl(|\partial_t^k \partial_{t,x} u^I|
+\frac{|\partial_t^k u^I|}{r}\Bigr) |\partial_t^j F^I|\:dx\:ds
+\sum_{j\le N} \|\partial_t^j G\|^2_{L^2_tL^2_x(S_t)}
\\+\sum_{j,k\le N} \int_0^t\int_\ext \Bigl|(\partial_\gamma h^{IJ,\gamma\delta}) \partial_t^j \partial_\delta u^J
\Bigr|\Bigl(|\partial_t^k \partial_{t,r} u^I|+\frac{|\partial_t^k u^I|}{r}\Bigr)\:dx\:ds
\\+\sum_{j,k\le N} \int_0^t \int_\ext \Bigl|(\partial_{t,r} h^{IJ,\gamma\delta}) \partial_t^j \partial_\gamma u^J
\partial_t^k \partial_\delta u^I\Bigr| 
\:dx\:ds
\\+\sum_{j,k\le N} \int_0^t \int_\ext \frac{|h|}{\langle x\rangle} |\partial_t^j u'|\Bigl(|\partial_t^k u'|
+\frac{|\partial_t^k u|}{r}\Bigr)\:dx\:ds + \sum_{|\mu|\le N-1} \|\partial^\mu \Box u\|_{L^2_tL^2_x(S_t)}.
\end{multline}
It should be noted that we now require additional smoothness of the boundary of $\mathcal{K}$, rather than 
$C^1$ as in Proposition \ref{prope.1} and Proposition \ref{prope.2}.

We will need a similar estimate involving the scaling vector field as well as derivatives.  In order to 
obtain this, we employ a technique from \cite{MS1} which introduces a modified scaling vector field
$\tilde{L}=t\partial_t + \eta(x)r\partial_r$ where $\eta$ is a smooth function with $\eta(x)\equiv 0$ for
$x\in\mathcal{K}$ and $\eta(x)\equiv 1$ for $|x|\ge 1$.  Here, of course, we are relying on the assumption
that $\mathcal{K}\subset \{|x|\le 1\}$.  

We will look to bound
\begin{equation}
\label{e.25}
\sum_{\substack{|\mu|+k\le N\\k\le K}} \|L^k \partial^\mu u'\|_{L^2_tL^2_x(S_t\cap \{|x| < 1\})}.
\end{equation}
By elliptic regularity, this is
\begin{equation}
\label{e.26} \begin{split}
&\lesssim \sum_{\substack{j+k\le N\\k\le K}} \|L^k \partial_t^j u'\|_{L^2_tL^2_x(S_t\cap \{|x|<3/2\})}
+\sum_{\substack{|\mu|+k\le N-1\\k\le K}} \|L^k \partial^\mu \Box u\|_{L^2_tL^2_x(S_t)}
\\&\lesssim \sum_{\substack{j+k\le N\\k\le K}} \|(\tilde{L}^k \partial_t^j u)'\|_{L^2_tL^2_x
(S_t\cap \{|x|<3/2\})} + \sum_{\substack{|\mu|+k\le N\\k\le K-1}} \|L^k \partial^\mu u'
\|_{L^2_tL^2_x(S_t\cap \{|x|<3/2\})} 
\\&\qquad\qquad\qquad\qquad\qquad\qquad\qquad\qquad\qquad\qquad
+ \sum_{\substack{|\mu|+k\le N-1\\k\le K}} \|L^k
\partial^\mu \Box u\|_{L^2_tL^2_x(S_t)}.
\end{split}
\end{equation}

If $P=P(t,x,D_t,D_x)$ is a differential operator, we fix the notation (as in \cite{MS1}):
$$[P,h^{\gamma\delta}\partial_\gamma\partial_\delta]u=\sum_{1\le I,J\le D} \sum_{0\le \gamma,\delta\le 3}
[P,h^{IJ,\gamma\delta}\partial_\gamma\partial_\delta]u^J.$$
Since
\begin{align*}
[\Box_h, \tilde{L}^k \partial_t^j]u&=[\Box, \tilde{L}^k \partial_t^j]u + [h^{\gamma\delta}\partial_\gamma
\partial_\delta, \tilde{L}^k \partial_t^j]u\\
&=[\Box, L^k]\partial_t^j u + [\Box,(\tilde{L}^k-L^k)]\partial_t^ju +[h^{\gamma\delta}\partial_\gamma
\partial_\delta,\tilde{L}^k \partial_t^j]u
\end{align*}
and since $\tilde{L}^k \partial^j_t u$ satisfies the Dirichlet
boundary condition, in order to bound the
first term in the right side of \eqref{e.26} we can apply
\eqref{e.23} with $F$ replaced by
$$\tilde{L}^k \partial_t^j F + [h^{\gamma\delta}\partial_\gamma\partial_\delta, \tilde{L}^k \partial_t^j]
u + [\Box, L^k] \partial_t^j u$$
and $G$ by
$$\tilde{L}^k \partial_t^j G + [\Box, \tilde{L}^k -L^k]\partial_t^j u,$$
which is supported in $|x|<2$.  Thus, it follows that
\begin{multline}
\label{e.27}
\sum_{\substack{|\mu|+k\le N\\k\le K}} \|L^k \partial^\mu u'\|^2_{L^2_tL^2_x(S_t\cap\{|x|<1\})}
\lesssim 
\sum_{\substack{|\mu|+k\le N\\k\le K}} \|L^k \partial^\mu u'(0,\cd)\|_2^2
\\+ \sum_{\substack{|\mu|+j\le N\\j\le K}} \sum_{\substack{|\nu|+k\le N\\k\le K}}
\int_0^t \int_\ext \Bigl(|L^j \partial^\mu \partial u^I|+\frac{|L^j \partial^\mu u^I|}{r}\Bigr)
|L^k \partial^\nu F^I|\:dx\:ds
\\+\sum_{\substack{|\mu|+j\le N\\j\le K}} \sum_{\substack{l+k\le N\\k\le K}}
\int_0^t \int_\ext \Bigl(|L^j \partial^\mu \partial u^I|+\frac{|L^j \partial^\mu u^I|}{r}\Bigr)
|[h^{IJ,\gamma\delta}\partial_\gamma\partial_\delta, \tilde{L}^k \partial_t^l] u^J|\:dx\:ds
\\+\sum_{\substack{|\mu|+j\le N\\j\le K}} \sum_{\substack{l+k\le N-1\\k\le K-1}}
\int_0^t \int_\ext \Bigl(|L^j \partial^\mu \partial u^I|+\frac{|L^j \partial^\mu u^I|}{r}\Bigr)
|L^k \partial_t^l (\Box u)^I|\:dx\:ds
\\+\sum_{\substack{|\mu|+j\le N\\j\le K}} \sum_{\substack{|\nu|+k\le N\\k\le K}}
\int_0^t \int_\ext \Bigl|(\partial_\gamma h^{IJ,\gamma\delta}) \partial_\delta(\tilde{L}^j \partial^\mu
u^J)\Bigr| \Bigl(|L^k \partial^\nu \partial_{t,r} u^I|+\frac{|L^k \partial^\nu u^I|}{r}\Bigr)\:dx\:ds
\\+\sum_{\substack{|\mu|+j\le N\\j\le K}}\sum_{\substack{|\nu|+k\le N\\k\le K}}
\int_0^t \int_\ext \Bigl|(\partial_{t,r} h^{IJ,\gamma\delta})\partial_\gamma(\tilde{L}^j \partial^\mu  u^J)
\partial_\delta(\tilde{L}^k \partial^\nu u^I)\Bigr|\:dx\:ds
\\+\sum_{\substack{|\mu|+j\le N\\j\le K}}\sum_{\substack{|\nu|+k\le N\\k\le K}}
\int_0^t \int_\ext \frac{|h|}{\langle x\rangle} |L^j \partial^\mu u'|\Bigl(|L^k \partial^\nu u'|
+\frac{|L^k \partial^\nu u|}{r}\Bigr)\:dx\:ds
\\+\sum_{\substack{|\mu|+k\le N\\k\le K}} \|L^k \partial^\mu G\|^2_{L^2_tL^2_x(S_t)}
\\+\sum_{\substack{|\mu|+k\le N-1\\k\le K}}\|L^k \partial^\mu \Box u\|^2_{L^2_tL^2_x(S_t)}
+\sum_{\substack{|\mu|+k\le N\\k\le K-1}} \|L^k \partial^\mu u'\|^2_{L^2_tL^2_x(S_t\cap \{|x|
<3/2\})},
\end{multline}
for solutions $u$ to \eqref{e.22}.
If we argue recursively, the same bound holds with the last term replaced by
$$\sum_{|\mu|\le N} \|\partial^\mu u'\|^2_{L^2_tL^2_x(S_t\cap \{|x|<2\})}.$$
Thus, in order to control this last term, we may apply \eqref{e.24}.  A similar argument can be used to 
bound
$$\sum_{|\mu|+k\le N} \|L^k \partial^\mu u'(t,\cd)\|_{L^2(\{|x|<1\})}.$$
Moreover, since $K\le N$ is
arbitrary, we have shown
\begin{lemma}\label{lemmae.3}
Suppose that $\mathcal{K}$ is a smooth, bounded, star-shaped obstacle as above.  Suppose
further that the perturbation terms $h^{IJ,\alpha\beta}$ are as above.  Then, if $u\in C^\infty$
solves \eqref{e.22} and vanishes for large $x$ for every $t$ and $G$ is supported in $|x|<2$,
\begin{multline}
\label{e.28}
\sum_{|\mu|+j\le N} \|L^j \partial^\mu u'(t,\cd)\|_{L^2(\{|x|<1\})}+
\sum_{|\mu|+j\le N} \|L^j \partial^\mu u'\|^2_{L^2_tL^2_x(S_t\cap\{|x|<1\})}
\\\lesssim 
\sum_{|\mu|+j\le N} \|L^j \partial^\mu u'(0,\cd)\|_2^2
\\+ \sum_{|\mu|+ j\le N} \sum_{|\nu|+k\le N}
\int_0^t \int_\ext \Bigl(|L^j \partial^\mu \partial u^I|+\frac{|L^j \partial^\mu u^I|}{r}\Bigr)
|L^k \partial^\nu F^I|\:dx\:ds
\\+\sum_{|\mu|+j\le N} \sum_{l+k\le N}
\int_0^t \int_\ext \Bigl(|L^j \partial^\mu \partial u^I|+\frac{|L^j \partial^\mu u^I|}{r}\Bigr)
|[h^{IJ,\gamma\delta}\partial_\gamma\partial_\delta, \tilde{L}^k \partial_t^l] u^J|\:dx\:ds
\\+\sum_{|\mu|+j\le N} \sum_{l+k\le N-1}
\int_0^t \int_\ext \Bigl(|L^j \partial^\mu \partial u^I|+\frac{|L^j \partial^\mu u^I|}{r}\Bigr)
|L^k \partial_t^l (\Box u)^I|\:dx\:ds
\\+\sum_{|\mu|+j\le N} \sum_{|\nu|+k \le N}
\int_0^t \int_\ext \Bigl|(\partial_\gamma h^{IJ,\gamma\delta}) \partial_\delta(\tilde{L}^j \partial^\mu
u^J)\Bigr| \Bigl(|L^k \partial^\nu \partial_{t,r} u^I|+\frac{|L^k \partial^\nu u^I|}{r}\Bigr)\:dx\:ds
\\+\sum_{|\mu|+j\le N}\sum_{|\nu|+k\le N}
\int_0^t \int_\ext \Bigl|(\partial_{t,r} h^{IJ,\gamma\delta})\partial_\gamma(\tilde{L}^j \partial^\mu  u^J)
\partial_\delta(\tilde{L}^k \partial^\nu u^I)\Bigr|\:dx\:ds
\\+\sum_{|\mu|+j\le N}\sum_{|\nu|+k\le N}
\int_0^t \int_\ext \frac{|h|}{\langle x\rangle} |L^j \partial^\mu u'|\Bigl(|L^k \partial^\nu u'|
+\frac{|L^k \partial^\nu u|}{r}\Bigr)\:dx\:ds
\\+\sum_{|\mu|+j\le N} \|L^j \partial^\mu G\|^2_{L^2_tL^2_x(S_t)}
+\sum_{|\mu|+j\le N-1}\|L^j \partial^\mu \Box u\|^2_{L^2_tL^2_x(S_t)}
\\+\sum_{|\mu|+j\le N-1} \|L^j \partial^\mu \Box u(t,\cd)\|^2_2
\end{multline}
for any $N\ge 0$ and for every $t>0$.
\end{lemma}

We use this to show that a version of \eqref{e.23} holds when $u$ is replaced by $\Gamma^\mu u$
where $\Gamma=\{L,\Omega,\partial\}$ is the set of ``admissible'' vector fields.  With $\eta$ as in the definition
of $\tilde{L}$, we set $\tilde{\Omega}_{ij}=\eta(x)\Omega_{ij}$, $1\le
i<j\le 3$ and $\tilde{\partial}_i=\eta(x)\partial_i$, $i=
1,2,3$.  Similarly, we set $\tilde{\Gamma}=\{\tilde{L},\tilde{\Omega}, \tilde{\partial}_i, \partial_t\}$ to
be the set of boundary-preserving vector fields.

\begin{theorem}
\label{theoreme.4}
Suppose $\mathcal{K}$ is a smooth, bounded, star-shaped obstacle as above.  Suppose further that
the perturbation terms $h^{IJ,\alpha\beta}$ are as above.  Then, if $u\in C^\infty$ solves \eqref{e.22}
and vanishes for large $x$ for every $t$ and $G$ is supported in $|x|<2$,
\begin{multline}
\label{e.29}
\sum_{|\mu|\le N} \|\langle x\rangle^{-1/2-} \Gamma^\mu u'\|^2_{L^2_tL^2_x(S_t)} + (\log(2+t))^{-1} 
\sum_{|\mu|\le N} \|\langle x\rangle^{-1/2} \Gamma^\mu u'\|^2_{L^2_tL^2_x(S_t)}
\\+\sum_{|\mu|\le N} \|\Gamma^\mu u'(t,\cd)\|_2
\lesssim \sum_{|\mu|\le N} \|\Gamma^\mu u'(0,\cd)\|^2_2
\\+\sum_{|\mu|,|\nu|\le N} \int_0^t \int_\ext \Bigl(|\Gamma^\mu \partial u^I|+\frac{|\Gamma^\mu u^I|}{r}
\Bigr) |\Gamma^\nu F^I|\:dx\:ds
\\+\sum_{|\mu|,|\nu|\le N} \int_0^t \int_\ext \Bigl(|\Gamma^\mu \partial u^I| + \frac{|\Gamma^\mu u^I|}{r}
\Bigr) |[h^{IJ,\gamma\delta}\partial_\gamma\partial_\delta, \Gamma^\nu]u^J|\:dx\:ds\\
\\+\sum_{|\mu|\le N,|\nu|\le N-1} \int_0^t \int_\ext \Bigl(|\Gamma^\mu \partial u^I| + 
\frac{|\Gamma^\mu u^I|}{r}\Bigr) |\Gamma^\nu (\Box u)^I|\:dx\:ds\\
+\sum_{|\mu|,|\nu| \le N} \int_0^t \int_\ext
|(\partial_\gamma h^{IJ,\gamma\delta}) \partial_\delta (\Gamma^\mu u^J)|
\Bigl(|\Gamma^\nu \partial u^I|+\frac{|\Gamma^\nu u^I|}{r}\Bigr)\:dx\:ds 
\\+\sum_{|\mu|,|\nu|\le N} \int_0^t \int_\ext \Bigl|(\partial_{t,r} h^{IJ,\gamma\delta}) \partial_\gamma
(\Gamma^\mu u^J) \partial_\delta (\Gamma^\nu u^I)\Bigr|\:dx\:ds
\\+\sum_{\substack{|\mu|+|\sigma|\le N\\|\nu|\le N}} \int_0^t 
\int_{\ext\cap \{|x|<1\}} |\Gamma^\sigma 
h| |\Gamma^\mu u'|
\Bigl(|\Gamma^\nu u'| + \frac{|\Gamma^\nu u|}{r}\Bigr)\:dx\:ds
\\+\sum_{|\mu|,|\nu|\le N} \int_0^t \int_\ext \frac{|h|}{\langle x\rangle} |\Gamma^\mu u'|\Bigl(
|\Gamma^\nu u'| + \frac{|\Gamma^\nu u|}{r}\Bigr)\:dx\:ds\\
\sum_{|\mu|\le N} \|\Gamma^\mu G\|^2_{L^2_tL^2_x(S_t)} + \sum_{|\mu|\le N-1} \|\Gamma^\mu \Box
u\|^2_{L^2_tL^2_x(S_t)} + \sum_{|\mu|\le N-1} \|\Gamma^\mu \Box u(t,\cd)\|^2_2
\end{multline}
for any fixed $N\ge 0$ and any $t\ge 0$.
\end{theorem}

To show this, we argue inductively in $N$, and the case $N=0$ clearly follows from \eqref{e.23}.  Let
us show the bound for the first term on the left side.  Similar arguments will yield the full estimate.

We begin by observing that
\begin{multline}
\label{e.30}
\sum_{|\mu|\le N} \|\langle x\rangle^{-1/2-} \Gamma^\mu u'\|_{L^2_tL^2_x(S_t)}
\lesssim \sum_{|\mu|\le N-1} \|\langle x\rangle^{-1/2-} \Gamma^\mu (\tilde{\Gamma}u)'\|_{L^2_tL^2_x(S_t)}
\\+ \sum_{\substack{|\mu|+j\le N}} \|L^j \partial^\mu u'\|_{L^2_tL^2_x(S_t\cap\{|x|<1\})}.
\end{multline}
The bound for the last term clearly follows from \eqref{e.28}.

To estimate the first term in the right side of \eqref{e.30}, we begin by noticing that
\begin{equation}
\label{e.31}\begin{split}
\Box_h \tilde{\Gamma}u &= \tilde{\Gamma}\Box_h u + [\Box_h, \tilde{\Gamma}]u\\
&=\tilde{\Gamma}\Box_h u + [\Box, \tilde{\Gamma}]u + [h^{\gamma\delta}\partial_\gamma\partial_\delta,
\tilde{\Gamma}] u\\
&=\tilde{\Gamma}\Box_h u + [\Box, \Gamma]u + [\Box, \tilde{\Gamma}-\Gamma]u + [h^{\gamma\delta}
\partial_\gamma\partial_\delta,\tilde{\Gamma}]u.
\end{split}
\end{equation}
Thus, we will apply the inductive hypothesis to $\tilde{\Gamma}u$ with $F$ replaced by
$$\tilde{\Gamma} F + [\Box,\Gamma]u + [h^{\gamma\delta}\partial_\gamma\partial_\delta, \tilde{\Gamma}]u$$
and $G$ by the compactly supported function
$$\tilde{\Gamma} G + [\Box, \tilde{\Gamma}-\Gamma]u.$$
It follows that the first term in the right side of \eqref{e.30} is dominated by the right side of
\eqref{e.29} plus
$$\sum_{|\mu|+j\le N} \|L^j \partial^\mu u'\|^2_{L^2_tL^2_x(S_t\cap\{|x|<1\})}$$
since the coefficients of $Z$ are $O(1)$ in $\{|x|<1\}$.  Thus, another application
of \eqref{e.28} yields the desired estimate.

\bigskip
\newsection{Decay estimates}

Classically, the necessary decay to prove long-time existence is afforded to us by the 
Klainerman-Sobolev inequalities (see \cite{K2}; see also \cite{H},\cite{S}).  These inequalities, however, require the
use of the Lorentz rotations which does not seem permissible in the current setting.  In order to
get around this, we will rely on decay in $|x|$ (which meshes well with the KSS estimates from the
previous section) obtained by a weighted Sobolev inequality and
decay in $t-|x|$ that follows from (variants of) estimates of Klainerman and Sideris \cite{KS}.

\subsection{Null form estimates}
We begin by providing the well-known decay that is obtained when employing the null condition.  The
proof that we present is essentially from \cite{Si3}.
\begin{lemma}
\label{lemmad.1}
Assume that the null condition, \eqref{i.3}, is satisfied.
Let $c_0=\min\{c_I/2\,:\, I=1,\dots,D\}$.  Then, for  $|x|\ge c_0t/2$,
\begin{multline}
\label{d.1}
|B^{KK,\alpha\beta}_{K,\gamma} \partial_\gamma u^K \partial_\alpha\partial_\beta v^K|
\\\lesssim \frac{1}{\langle t+|x|\rangle} \Bigl[|\Gamma u^K||\partial^2 v^K| + |\partial u^K|
|\partial \Gamma v^K| + \langle c_K t-r\rangle |\partial u^K| |\partial^2 v^K|\Bigr],
\end{multline}
\begin{multline}\label{d.2}
|B^{KK,\alpha\beta}_{K,\gamma} \partial_\alpha\partial_\gamma u^K \partial_\beta v^K|
\\\lesssim \frac{1}{\langle t+|x|\rangle} \Bigl[|\Gamma v^K||\partial^2 u^K| + |\partial v^K|
|\partial \Gamma u^K| + \langle c_K t-r\rangle |\partial v^K| |\partial^2 u^K|\Bigr],
\end{multline}
and
\begin{multline}
\label{d.3}
|B^{KK,\alpha\beta}_{K,\gamma}\partial_\alpha u^K \partial_\beta v^K \partial_\gamma w^K|
\\\lesssim \frac{1}{\langle t+|x|\rangle} \Bigl[|\Gamma u^K||\partial v^K||\partial w^K|
+ |\partial u^K||\Gamma v^K||\partial w^K| + |\partial u^K||\partial v^K||\Gamma w^K|
\\+ \langle c_K t-r\rangle |\partial u^K||\partial v^K||\partial w^K|\Bigr].
\end{multline}
\end{lemma}

While \eqref{d.2} did not appear explicitly in \cite{Si3}, it is used there.  It is easily seen to
follow from the same argument as the other bounds.  We will show \eqref{d.1}.  The other estimates result
from similar arguments.

It suffices to consider the case $|(t,x)|\ge 1$ as the bounds are otherwise trivial.  We first decompose
the spatial gradient into radial and angular parts:
$$\nabla_x = \frac{x}{r}\partial_r - \frac{x}{r^2}\wedge \Omega,$$
where $\Omega=x\wedge \partial_x$ and $\wedge$ denotes the usual vector cross product.  Introducing
the operators $D^{\pm}=\frac{1}{2} (\partial_t \pm c_K \partial_r)$ and the null vectors $Y^{\pm}=(1,\pm x/c_Kr)$,
we can further decompose
$$(\partial_t,\nabla_x)=(Y^-D^-+Y^+D^+)-\Bigl(0,\frac{x}{r^2}\wedge\Omega\Bigr),$$
or alternately
\begin{align*}
\partial_{t,x}&=Y^-D^- - \frac{c_Kt-r}{c_Kt+r} Y^+D^- + \frac{c_K}{c_Kt+r} Y^+ L - 
\Bigl(0,\frac{x}{r^2}\wedge\Omega\Bigr)\\
&= Y^-D^-+R,\end{align*}
where 
\begin{equation}\label{null.1}
|Ru|\lesssim \langle r\rangle^{-1} |\Gamma u| + \frac{\langle c_Kt-r\rangle}{t+r} |\partial u|.
\end{equation}
Thus, we have
\begin{multline*}
B^{KK,\alpha\beta}_{K,\gamma} \partial_\gamma u^K \partial_\alpha\partial_\beta v^K
= B^{KK,\alpha\beta}_{K,\gamma} [Y^-_\alpha Y^-_\beta Y^-_\gamma D^- u^K (D^-)^2 v^K
+ R_\gamma u^K \partial_\alpha \partial_\beta v^K 
\\+ Y^-_\gamma D^- u^K R_\alpha \partial_\beta v^K
+Y^-_\gamma D^- u^K Y^-_\alpha D^- R_\beta v^K].
\end{multline*}
Since $(Y^-_0)^2/c_K^2 - (Y_1^-)^2-(Y_2^-)^2-(Y_3^-)^2=0$, by \eqref{i.3}, the first term must vanish.  The
remaining bounds follow from \eqref{null.1}.

\subsection{Weighted Sobolev estimates}

The first estimate is a now standard weighted Sobolev inequality.  See \cite{K2}.  The reader is also
encouraged to see \cite{KSS2} for the first example of how this decay can be paired with KSS estimates
to yield long time existence for nonlinear equations.

\begin{lemma}
\label{lemmad.2}
Suppose $h\in C^\infty(\R^3)$.  Then, for $R>1$
\begin{equation}
\label{d.4}
\|h\|_{L^\infty(R/2<|x|<R)} \lesssim R^{-1}\sum_{|\alpha|\le 2} \|Z^\alpha h\|_{L^2(R/4<|x|<2R)}.
\end{equation}
\end{lemma}

For $|x|\in (R/2,R)$, we apply Sobolev's estimate for $\R\times S^2$ to see that
$$|h(x)|\lesssim \sum_{|\alpha|+j\le 2} \Bigl(\int_{|x|-1/4}^{|x|+1/4} \int_{S^2} |\partial_r^j
\Omega^\alpha h(r\omega)|^2\:dr\:d\omega\Bigr)^{1/2}.$$
Since the volume element in $\R^3$ is a constant times $r^2\,dr\,d\omega$, this is dominated
by the right side of \eqref{d.4} as desired.

The second of the necessary Sobolev type estimates follows essentially from that in 
\cite{Si} (Lemma 3.3).
\begin{lemma}
Let $u\in C^\infty(\ext)$ and suppose that $u$ vanishes on $\bdy$ and for large $x$ for every $t$.  Then,
\label{lemmad.3}
\begin{equation}
\label{d.5}
r^{1/2}|u(t,x)| \lesssim \sum_{|\mu|\le 1} \|Z^\mu u'(t,\cd)\|_2.
\end{equation}
Moreover, 
\begin{equation}
\label{d.5'}
r^{1/2}\sum_{|\nu|\le N} |\Gamma^\nu u(t,x)| \lesssim \sum_{|\nu|\le N+1} \|\Gamma^\nu u'(t,\cd)\|_2
\end{equation}
for any $N\ge 0$.
\end{lemma}


We first note that the Dirichlet boundary condition allows us to control $u$ locally by $u'$.  Thus,
over $|x|\le 1$, the result follows trivially from the standard Sobolev estimates.

In the remaining region, $|x|\ge 1$, \eqref{d.5} is a consequence of the arguments in \cite{Si}.
We write $x=r\omega$ where $\omega\in S^2$ (and $d\omega$ denotes the surface measure of this
unit sphere).  We begin by noting that
$$r^2 |u(t,x)|^4 \lesssim \sum_{|\mu|\le 1} r^2 \|\Omega^\mu u(t,r\cd)\|^4_{L^4(S^2)}$$
follows from a basic Sobolev estimate.  By the fundamental theorem of calculus, 
H\"older's inequality, and the standard Sobolev estimate $\|h\|_6\lesssim \|\nabla h\|_2$, it follows
that
\begin{align*}
r^2 \int_{S^2} |v(t,r\omega)|^4\:d\omega &\lesssim r^2\int_r^\infty \int_{S^2} |\partial_r v(\rho\omega)|
|v(\rho\omega)|^3\:d\rho\:d\omega\\
&\lesssim \|\partial_r v(t,\cd)\|_2 \|v(t,\cd)\|_6^3 \lesssim \|\nabla_x v(t,\cd)\|_2^4
\end{align*}
which yields \eqref{d.5} when $v$ is replaced by $\Omega^\mu u$.

When $|x|\ge 1$, \eqref{d.5'} follows from the same argument as that for \eqref{d.5}.  Since the coefficients
of $Z$ are $O(1)$ for $|x|\le 1$, it only remains to show that
\begin{equation}\label{d.5''}
\sum_{|\nu|\le N} |L^\nu u(t,x)|\lesssim \sum_{|\nu|\le N+1} \|\Gamma^\nu u'(t,\cd)\|_2,
\quad |x|\le 1.\end{equation}
By Sobolev's lemma, we have that for $|x|\le 1$,
$$\sum_{|\nu|\le N} |L^\nu u(t,x)|\lesssim \sum_{|\nu|\le N+1} \|\Gamma^\nu u'(t,\cd)\|_2
+\sum_{|\nu|\le N} \|L^\nu u(t,\cd)\|_{L^2(\{|x|<2\})}.$$
The last term in the right side is
$$\lesssim \sum_{|\nu|\le N} \|(t\partial_t)^\nu u(t,\cd)\|_{L^2(\{|x|<2\})} + 
\sum_{|\nu|\le N-1} \|L^\nu u'(t,\cd)\|_{L^2(\{|x|<2\})}.$$
Since $\partial_t$ preserves the Dirichlet boundary conditions, it follows from the Fundamental Theorem
of Calculus that the former term is
$$\lesssim \sum_{|\nu|\le N} \|(t\partial_t)^\nu u'(t,\cd)\|_{L^2(\{|x|<2\})}
\lesssim \sum_{|\nu|+|\mu|\le N} \|L^\nu \partial^\mu u'(t,\cd)\|_{L^2(\{|x|<2\})}$$
which completes the proof of \eqref{d.5'}.

\subsection{Klainerman-Sideris estimates}
We finally present some estimates from \cite{KS} and some consequences of these estimates.  These estimates
are the ones that provide any required decay in the time variable $t$.

We begin with the following basic estimate from \cite{KS},
\begin{equation}
\label{d.6}
\langle c_Kt-r\rangle\Bigl(|\partial_t \partial u^K| + |\Delta u^K|\Bigr)\lesssim \sum_{|\mu|\le 1} |\Gamma^\mu
u'| + \langle t+r\rangle |\Box u|.
\end{equation}
Moreover, using integration by parts, it was shown that
\begin{equation}
\label{d.7}
\|\langle c_K t-r\rangle \partial^2 v^K(t,\cd)\|_2 \lesssim \sum_{|\mu|\le 1} \|\Gamma^\mu v'(t,\cd)\|_2
+ \|\langle t+r\rangle \Box v(t,\cd)\|_2
\end{equation}
when there is no boundary.  Moreover, if one applies the boundaryless analog of
\eqref{d.5} to $\langle c_Kt-r\rangle \partial v^K$ and
uses \eqref{d.7}, the following is obtained,
\begin{equation}
\label{d.7'}
r^{1/2}\langle c_Kt-r\rangle |\partial v^K(t,x)|\lesssim \sum_{|\mu|\le 2} \|\Gamma^\mu v'(t,\cd)\|_2
+\sum_{|\mu|\le 1} \|\langle t+r\rangle \Gamma^\mu \Box v(t,\cd)\|_2,
\end{equation}
which first appeared in Hidano and Yokoyama \cite{HidY}.

When there is a boundary, the integration by parts argument in \cite{KS} does not yield \eqref{d.7}.  We will,
however, require analogous estimates.  The first three are from \cite{MNS1}.  The first follows from applying
\eqref{d.7} to $\eta(x)u(t,x)$ where $\eta$ is a smooth cutoff that vanishes for $|x|\le 1$ and is identically
one when $|x|\ge 3/2$.  
\begin{multline}
\label{d.8}
\sum_{|\mu|\le N} \|\langle c_K t-r\rangle \Gamma^\mu \partial^2 u^K(t,\cd)\|_2
\lesssim \sum_{|\mu|\le N+1} \|\Gamma^\mu u'(t,\cd)\|_2 \\+ \sum_{|\mu|\le N} \|\langle t+r\rangle
\Gamma^\mu \Box u(t,\cd)\|_2 + t\sum_{|\mu|\le N} \|\Gamma^\mu u'(t,\cd)\|_{L^2(\{|x|<1\})}.
\end{multline}
Moreover, by combining \eqref{d.4}, \eqref{d.8}, and elliptic
regularity (cf. \cite{MNS1}), one obtains
\begin{multline}
\label{d.9}
r \langle c_K t-r\rangle \sum_{|\mu|\le N} |\Gamma^\mu \partial^2 u^K| \lesssim 
\sum_{|\mu|\le N+3} \|\Gamma^\mu u'(t,\cd)\|_2 \\+ \sum_{|\mu|\le N+2} \|\langle t+r \rangle
\Gamma^\mu \Box u(t,\cd)\|_2 + t \sum_{|\mu|\le N} \|\Gamma^\mu u'(t,\cd)\|_{L^2(\{|x|<1\})}.
\end{multline}
Finally, by applying \eqref{d.5} to the cutoff solution and using \eqref{d.8}, we can obtain the following
analog of the estimate from \cite{HidY}.
\begin{multline}
\label{d.10}
r^{1/2} \langle c_kt-r\rangle \sum_{|\mu|\le N} |\Gamma^\mu \partial u^K|
\lesssim \sum_{|\mu|\le N+2} \|\Gamma^\mu u'(t,\cd)\|_2 + \sum_{|\mu|\le N+1} \|\langle
t+r\rangle \Gamma^\mu \Box u(t,\cd)\|_2 
\\+t\sum_{|\mu|\le N} \|\Gamma^\mu u'(t,\cd)\|_{L^2(|x|<1)}.
\end{multline}

As in \cite{MS1}, in a region $|x|\ge (c_0/2)t$, the boundary terms are no longer required.
In particular, we have
\begin{multline}
\label{d.11}
\sum_{|\mu|\le N} \|\langle c_K t-r\rangle \partial^2 \Gamma^\mu u^K(t,\cd)\|_{L^2(|x|\ge c_0t/2)} 
\\\lesssim \sum_{|\mu|\le N+1} \|\Gamma^\mu u'(t,\cd)\|_2 + \sum_{|\mu|\le N} \|\langle t+r\rangle
\Gamma^\mu \Box u(t,\cd)\|_{L^2(|x|\ge c_0t/4)},
\end{multline}
\begin{multline}
\label{d.12}
r\langle c_K t-r\rangle \sum_{|\mu|\le N} |\partial^2 \Gamma^\mu u^K(t,x)|
\\\lesssim \sum_{|\mu|\le N+3} \|\Gamma^\mu u'(t,\cd)\|_2 + \sum_{|\mu|\le N+2} \|\langle t+r\rangle
\Gamma^\mu \Box u(t,\cd)\|_{L^2(|x|\ge c_0t/4)},
\end{multline}
and
\begin{multline}
\label{d.13}
r^{1/2} \langle c_K t-r\rangle \sum_{|\mu|\le N} |\partial \Gamma^\mu u^K(t,x)|
\\\lesssim \sum_{|\mu|\le N+2} \|\Gamma^\mu u'(t,\cd)\| + \sum_{|\mu|\le N+1} \|\langle t+r\rangle
\Gamma^\mu \Box u(t,\cd)\|_{L^2(|x|\ge c_0t/4)}.
\end{multline}
Indeed, we now fix $\eta\in C^\infty(\R^3)$ satisfying $\eta(x)\equiv 1$, $|x|>1/2$ and $\eta(x)\equiv 0$ for
$|x|<1/4$.  We then set $v(t,x)=\eta(x/(c_0\langle t\rangle))u(t,x)$ and apply \eqref{d.7} and \eqref{d.7'}.

\bigskip
\newsection{Global existence}

In this section, we prove our main result, Theorem \ref{theoremi.1}.  Here, we shall choose $N=30$, but this
is far from optimal.
As in \cite{Si3}, the proof proceeds by examining a coupling between a low-order energy and a higher-order energy.
\begin{equation}
\label{g.1}
\sum_{|\mu|\le 20} \Bigl(\|\Gamma^\mu u'(t,\cd)\|_2 + \|\langle x\rangle^{-5/8} \Gamma^\mu
u'\|_{L^2_tL^2_x(S_t)} 
\Bigr)
\le A\varepsilon
\end{equation}
\begin{equation}
\label{g.2}
\sum_{|\mu|\le 30} \Bigl(\|\Gamma^\mu u'(t,\cd)\| +  \|\langle x\rangle^{-5/8}
\Gamma^\mu u'\|_{L^2_tL^2_x(S_t)} 
\Bigr)
\le B\varepsilon (1+t)^{c\varepsilon}.
\end{equation}
Here, $A$ is chosen to be 10 times greater than the square root of the implicit constant in \eqref{e.29}.  The
exponent $5/8$ was chosen to make the argument explicit.  The same argument would hold for sufficiently small $\varepsilon$ for
any exponent $p$ with $1/2<p<3/4$.

There are two steps required in order to complete the continuity argument:
\begin{enumerate}
\item[$(i.)$] Show \eqref{g.1} holds with $A$ replaced by $A/2$,
\item[$(ii.)$] Show that \eqref{g.2} follows from \eqref{g.1}.
\end{enumerate}

Throughout the remainder of the argument, we will be applying \eqref{e.29} with 
$$h^{IJ,\alpha\beta} = - B^{IJ,\alpha\beta}_{K,\gamma} \partial_\gamma u^K$$
and $F=G=0$.

\subsection{Preliminaries} 
Before beginning the proofs of $(i.)$ and $(ii.)$, we establish some 
preliminary estimates.  These are shown assuming \eqref{g.1}, and both are used
to control terms that appear after applications of the decay estimates.

The first is a lower order version.  We will establish:
\begin{equation}
\label{g.3}
\sum_{|\mu|\le 19} \|\langle t+r\rangle \Gamma^\mu  \Box u(t,\cd)\|_2 \lesssim \varepsilon^2
+ t \varepsilon \sum_{|\mu|\le 11} \|\Gamma^\mu
u'(t,\cd)\|_{L^2(|x|<1)}.
\end{equation}

The left side of \eqref{g.3} is clearly controlled by
$$\sum_{|\mu|\le 11,|\nu|\le 20} \|\langle t+r \rangle \Gamma^\mu u'(t,\cd) \Gamma^\nu u'(t,\cd)\|_2.$$
When the norm is taken over $|x|\ge c_0t/2$, we can apply \eqref{d.4} and \eqref{g.1} to see that this is $O(\varepsilon^2)$.
When the norm is over $|x|\le c_0t/2$, we apply \eqref{d.10} to see that this is
\begin{multline*}
\sum_{|\nu|\le 20} \|\Gamma^\nu u'(t,\cd)\|_2 \Bigl(
\sum_{|\mu|\le 14} \|\Gamma^\mu u'(t,\cd)\|_2 \\+ \sum_{|\mu|\le 13} \|\langle t+r\rangle \Gamma^\mu
\Box u(t,\cd)\|_2 + t \sum_{|\mu|\le 11} \|\Gamma^\mu u'(t,\cd)\|_{L^2(|x|<1)}\Bigr)
\\\lesssim \varepsilon^2 + \varepsilon \sum_{|\mu|\le 13} \|\langle t+r\rangle \Gamma^\mu \Box u(t,\cd)\|_2
+\varepsilon t\sum_{|\mu|\le 11} \|\Gamma^\mu
u'(t,\cd)\|_{L^2(|x|<1)}
.
\end{multline*}
The last inequality follows from \eqref{g.1}.
Since the second term on the right can be bootstrapped if $\varepsilon$ is sufficiently small, we see that
this yields \eqref{g.3}.

From this proof, it is easy to see that we also have
\begin{equation}
\label{g.3'}
\sum_{|\mu|\le 19} \|\langle t+r\rangle \Gamma^\mu \Box u(t,\cd)\|_{L^2(|x|\ge c_0t/4)} \lesssim \varepsilon^2.
\end{equation}

We will additionally require the related higher order estimate
\begin{equation}
\label{g.4}
\sum_{|\mu|\le 29} \|\langle t+r\rangle \Gamma^\mu \Box u(t,\cd)\|_2
\lesssim \Bigl(\varepsilon  +  t \sum_{|\mu|\le 15}
\|\Gamma^\mu u'(t,\cd)\|_{L^2(|x|<1)}\Bigr)\sum_{|\nu|\le 30} \|\Gamma^\nu u'(t,\cd)\|_2.
\end{equation}

Plugging in our nonlinearity in the left, this is
$$\lesssim \sum_{|\mu|\le 30, |\nu|\le 15} \|\langle t+r\rangle \Gamma^\mu u' \Gamma^\nu u'\|_2.
$$
If $|x|\ge c_0t/2$, applying \eqref{d.4} and \eqref{g.1} results in the bound
$$\lesssim \varepsilon \sum_{|\mu|\le 30} \|\Gamma^\mu u'(t,\cd)\|_2.$$
If $|x|<c_0t/2$, we apply \eqref{d.10} to see that this is
\begin{multline*}
\lesssim \sum_{|\mu|\le 17} \|\Gamma^\mu u'(t,\cd)\| \sum_{|\nu|\le 30} \|
\Gamma^\nu u'(t,\cd)\|_2 + \sum_{|\mu|\le 16} \|\langle t+r\rangle \Gamma^\mu \Box u(t,\cd)\|_2
\sum_{|\nu|\le 30} \|\Gamma^\nu u'(t,\cd)\|_2
\\+ \sum_{|\mu|\le 30} \|\Gamma^\mu u'(t,\cd)\|_2
\sum_{|\nu|\le 15} t\|\Gamma^\nu u'(t,\cd)\|_{L^2(|x|<1)}.\end{multline*}
By applying \eqref{g.1} and \eqref{g.3}, we indeed see that \eqref{g.4} follows.

Again, the same proof also yields
\begin{equation}
\label{g.4'}
\sum_{|\mu|\le 29} \|\langle t+r\rangle \Gamma^\mu \Box u(t,\cd)\|_{L^2(|x|\ge c_0t/2)}
\lesssim \varepsilon \sum_{|\mu|\le 30} \|\Gamma^\mu u'(t,\cd)\|_2.
\end{equation}

\subsection{Low order energy}
Here, while assuming \eqref{g.1} and \eqref{g.2}, we must show that \eqref{g.1} holds with
$A$ replaced by $A/2$.  Using \eqref{i.4} and \eqref{e.29}, the square of the left side of \eqref{g.1} is
easily seen to be
\begin{multline}
\label{l.1}
\le (A\varepsilon/10)^2 + C\sum_{|\mu|\le 20} \sum_{|\nu|\le 19, |\sigma|\le 10} 
\int_0^t \int_\ext |\Gamma^\mu \partial u^I| |\tilde{B}^{IJ,\alpha\beta}_{K,\gamma} \partial_\gamma 
(\Gamma^\sigma u^K) \partial_\alpha\partial_\beta (\Gamma^\nu u^J)|\:dx\:ds\\
+C \sum_{|\mu|,|\sigma|\le 20} \sum_{|\nu|\le 10} \int_0^t \int_\ext |\Gamma^\nu \partial u^I|
|\tilde{B}^{IJ,\alpha\beta}_{K,\gamma}
\partial_\gamma (\Gamma^\sigma u^K) \partial_\alpha\partial_\beta (\Gamma^\nu u^J)|\:dx\:ds
\\+C \sum_{|\mu|,|\nu|\le 20}  \int_0^t \int_\ext |\Gamma^\mu \partial u^I|
|\tilde{B}^{IJ,\alpha\beta}_{K,\gamma} \partial_\alpha\partial_\gamma u^K \partial_\beta (\Gamma^\nu u^J)|
\:dx\:ds
\\+C \sum_{|\mu|,|\nu|\le 20} \int_0^t \int_\ext |\tilde{B}^{IJ,\alpha\beta}_{K,\gamma} \partial_{\gamma}
\partial_{t,r} u^K \partial_\alpha (\Gamma^\mu u^I) \partial_\beta (\Gamma^\nu u^J)|\:dx\:ds
\\+C \sum_{|\mu|\le 20} \sum_{|\nu|\le 19}\sum_{|\sigma|\le 11} \int_0^t \int_\ext
\frac{|\Gamma^\mu u|}{r} |\Gamma^\nu u''| |\Gamma^\sigma u'|\:dx\:ds
\\+C\sum_{|\mu|,|\nu|\le 20} \sum_{|\sigma|\le 11} \int_0^t \int_\ext
\frac{|\Gamma^\mu u|}{r} |\Gamma^\nu u'||\Gamma^\sigma u''|\:dx\:ds
\\+C \sum_{|\mu|,|\nu|\le 20} \sum_{|\sigma|\le 11} \int_0^t \int_{\ext\cap \{|x|<1\}}
|\Gamma^\sigma u'| |\Gamma^\mu u'|\Bigl(|\Gamma^\nu u'| + \frac{|\Gamma^\nu u|}{r}\Bigr)\:dx\:ds
\\+C \sum_{|\mu|,|\nu|\le 20} \int_0^t \int_\ext \frac{1}{\langle x\rangle} |u'||\Gamma^\mu u'| 
\Bigl(|\Gamma^\nu u'|+\frac{|\Gamma^\nu u|}{r}\Bigr)\:dx\:ds\\
+C\sum_{|\mu|\le 20}\sum_{|\nu|\le 11} \||\Gamma^\nu u'|\,|\Gamma^\mu u'|\|^2_{L^2_tL^2_x(S_t)}
+C\sum_{|\mu\le 20}\sum_{|\nu|\le 11} \||\Gamma^\nu u'(t,\cd)| \,|\Gamma^\mu u'(t,\cd)|\|^2_2.
\end{multline}
Due to constants that are introduced when $\Gamma^\mu$ commutes with $\partial_\alpha$, the 
coefficients $B^{IJ,\alpha\beta}_{K,\gamma}$ become new constants $\tilde{B}^{IJ,\alpha\beta}_{K,\gamma}$.
By, e.g.,  Lemma 4.1 of \cite{Si3}, it is known that $\Gamma$ preserves the null forms.  Thus, if $B^{KK,\alpha\beta}
_{K,\gamma}$ satisfies \eqref{i.3}, then so do the $\tilde{B}^{KK,\alpha\beta}_{K,\gamma}$.

In order to complete the proof, we will show that every term in \eqref{l.1} except for the first
is $O(\varepsilon^3)$ if $\varepsilon$ is sufficiently small.  By \eqref{d.4}, a Hardy 
inequality, and the Schwarz inequality, the eighth and ninth terms above are dominated by
$$C\sum_{|\mu|\le 20} \Bigl(\sup_{0\le s\le t} \|\Gamma^\mu u'(s,\cd)\|_2\Bigr) \sum_{|\nu|\le 20}
\|\langle x\rangle^{-1} \Gamma^\nu u'\|^2_{L^2_tL^2_x(S_t)}$$
which is $O(\varepsilon^3)$ by \eqref{g.1}.  We may similarly apply \eqref{d.4} to see that the last two
terms are
$$\lesssim \sum_{|\mu|\le 20}\Bigl(\sup_{0\le s\le t} \|\Gamma^\mu u'(s,\cd)\|^2_2\Bigr)
\Bigl(\sum_{|\nu|\le 20} \|\langle x\rangle^{-1} \Gamma^\nu u'\|^2_{L^2_tL^2_x(S_t)}
+\sum_{|\nu|\le 20} \|\Gamma^\nu u'(t,\cd)\|_2^2\Bigr)$$
which is $O(\varepsilon^4)$ by \eqref{g.1}.

It thus suffices to show that the second through the seventh terms in \eqref{l.1} are $O(\varepsilon^3)$.
We shall examine
the regions $|x|\le c_0 s/2$ and $|x|\ge c_0 s/2$ separately.  Here, $c_0$ is as in Lemma \ref{lemmad.1}.

\subsubsection{In the region $|x|\le c_0 s/2$:}  This is the easier case.  We look at the remaining
terms in \eqref{l.1} when the spatial integrals are restricted to $|x|\le c_0s/2$.  By the Schwarz inequality
and a Hardy inequality, these terms are
\begin{equation}
\label{l.2}\begin{split}
&\lesssim \Bigl(\sum_{|\mu|\le 20} \sup_{0\le s\le t} \|\Gamma^\mu u'(s,\cd)\|_2\Bigr) \Bigl[\int_0^t
 \sum_{|\nu|\le 20, |\sigma|\le 10} \|\Gamma^\sigma u''\, \Gamma^\nu u'(s,\cd)\|_{L^2(|x|\le c_0s/2)}\:ds
\\&\qquad\qquad\qquad\qquad\qquad\qquad\qquad
+\int_0^t \sum_{|\nu|\le 19, |\sigma|\le 10} \|\Gamma^\sigma u'\, \Gamma^\nu u''(s,\cd)\|_{L^2(|x|\le c_0s/2)}\:ds\Bigr]\\
&\lesssim \varepsilon \Bigl[\int_0^t 
 \sum_{|\nu|\le 20, |\sigma|\le 10} \|\Gamma^\sigma u''\, \Gamma^\nu u'(s,\cd)\|_{L^2(|x|\le c_0s/2)}\:ds
\\&\qquad\qquad\qquad\qquad\qquad\qquad\qquad
+\int_0^t\sum_{|\nu|\le 19, |\sigma|\le 10} \|\Gamma^\sigma u'\, \Gamma^\nu u''(s,\cd)\|_{L^2(|x|\le c_0s/2)}\:ds\Bigr].
\end{split}
\end{equation}
The second inequality follows from \eqref{g.1}.

For the first term on the right of \eqref{l.2}, we apply \eqref{d.9} and \eqref{g.3} to see that this
term is
\begin{multline*}
\lesssim \varepsilon^2 \sum_{|\mu|\le 20} \int_0^t \frac{1}{\langle s\rangle} \|\langle x\rangle^{-1}
\Gamma^\mu u'(s,\cd)\|_2\:ds \\+ \varepsilon \sum_{|\mu|\le 20, |\nu|\le 11} 
\int_0^t \|\Gamma^\nu u'(s,\cd)\|_{L^2(|x|<1)}
\|\langle x\rangle^{-1} \Gamma^\mu u'(s,\cd)\|_2\:ds. 
\end{multline*}
Here, we have also applied the bound \eqref{g.1}.  By the Schwarz inequality and \eqref{g.1}, this is indeed 
$O(\varepsilon^3)$.
For the second term on the right of \eqref{l.2}, we apply \eqref{d.4} and \eqref{d.8} to the two factors
respectively, yielding the same bound as above
for the first term in the right of \eqref{l.2}, but with the weights $\langle x\rangle^{-1}$ 
replaced by $\langle x\rangle^{-1+}$.  
Thus, this
term is also $O(\varepsilon^3)$ as desired.

\subsubsection{In the region $|x|> c_0 s/2$:}
We first show that the sixth and seventh terms in \eqref{l.1} are $O(\varepsilon^3)$ when the 
spatial integral is taken over $|x|>c_0s/2$.  Indeed, we can apply \eqref{d.5'} to see that these terms
are
$$
\lesssim \sum_{|\mu|\le 21} \sum_{|\nu|\le 20} \int_0^t \frac{1}{(1+s)^{1/4}} \|\Gamma^\mu u'(s,\cd)\|_2
\|\langle x\rangle^{-5/8} \Gamma^\nu u'(s,\cd)\|^2_2\:ds.
$$
Thus, by \eqref{g.1} and \eqref{g.2}, this is indeed
$O(\varepsilon^3)$ for $\varepsilon$ sufficiently small.

For the remaining terms (the second, third, fourth, and fifth terms in \eqref{l.1}), there are two cases: $(1)$ when all
three wave speeds are the same, $(I,J,K)=(I,I,I)$, and $(2)$ when there is a wave speed that is distinct from the other
two.

In case $(1)$, the null form bounds \eqref{d.1}, \eqref{d.2}, and \eqref{d.3} apply.  In the region $|x|\ge c_0s/2$,
these terms are 
\begin{multline}
\label{l.3}
\lesssim \int_0^t \int_{|x|\ge c_0s/2} \frac{1}{\langle s+r\rangle} \Bigl(\sum_{|\mu|\le 11} |\Gamma^\mu u'|
\sum_{|\nu|\le 20} |\Gamma^\nu u| \sum_{|\sigma|\le 20} |\Gamma^\sigma u'|\Bigr)\:dx\:ds\\
+\int_0^t \int_{|x|\ge c_0s/2} \frac{\langle c_K s-r\rangle}{\langle s+r\rangle} \Bigl(\sum_{|\mu|\le 11} |\Gamma^\mu 
\partial u^K|\Bigr)^2 \sum_{|\nu|\le 20} |\Gamma^\nu \partial u^K|\:dx\:ds.
\end{multline}
Applying \eqref{d.5'}, it follows as above that the first term is
$$\lesssim \int_0^t \int_{|x|\ge c_0s/2} \frac{1}{(1+s)^{1/4}} \sum_{|\mu|\le 21} \|\Gamma^\mu u'(s,\cd)\|_2
\sum_{|\mu|\le 20} \|\langle x\rangle^{-5/8} \Gamma^\mu u'(s,\cd)\|^2_2\:ds$$
which is easily seen to be $O(\varepsilon^3)$ using \eqref{g.1} and \eqref{g.2}.  For the second term
in \eqref{l.3}, we apply \eqref{d.13}, \eqref{g.1}, and \eqref{g.3'} to see that it is
$$
\lesssim \varepsilon \int_0^t \int_{|x|\ge c_0s/2} \frac{1}{\langle x\rangle^{-3/2}} 
\sum_{|\mu|\le 20} |\Gamma^\mu u'(s,x)|^2
\:dx\:ds,
$$
which is $O(\varepsilon^3)$ by \eqref{g.1}.  This concludes the proof of the same speed case $(1)$.

We next examine case $(2)$, the multiple speed case.  Here, we must bound
\begin{equation}
\label{l.4}
\int_0^t \int_{|x|\ge c_0s/2} \sum_{|\mu|\le 20} |\Gamma^\mu \partial u^I|\sum_{|\nu|\le 20} |\Gamma^\nu \partial u^J|
\sum_{|\sigma|\le 20} |\Gamma^\sigma \partial u^K|\:dx\:ds
\end{equation}
with $(I,K)\neq (K,J)$.  For simplicity of exposition, we assume $I\neq K$, $I=J$.  The other cases follow from
symmetric arguments.  We fix $\delta<|c_I-c_K|/3$.  Thus, $\{|x|\in [(c_I-\delta)s,(c_I+\delta)s]\}\cap
\{|x|\in [(c_K-\delta)s,(c_K+\delta)s]\}=\emptyset$, and it suffices to show that \eqref{l.4} is $O(\varepsilon^3)$
when the spatial integral is taken over the complements of these sets separately.  We will show the bound
over $\{|x|\not\in [(c_K-\delta)s,(c_K+\delta)s]\}\cap \{|x|\ge c_0s/2\}$.  Again, the remainder of the necessary
argument follows symmetrically.

If we apply \eqref{d.13} and \eqref{g.4'}, we have
\begin{multline*}
\int_0^t \int_{\{|x|\not\in [(c_K-\delta)s,(c_K+\delta)s]\}\cap \{|x|\ge c_0s/2\}}
\Bigl(\sum_{|\mu|\le 20} |\Gamma^\mu \partial u^I|\Bigr)^2 \sum_{|\nu|\le 20} |\Gamma^\nu \partial u^K|\:dx\:ds
\\\lesssim \int_0^t \frac{1}{(1+s)^{1/4}} \sum_{|\nu|\le 30} \|\Gamma^\nu u'(s,\cd)\|_2
\int_\ext \Bigl(\langle x\rangle^{-5/8} \sum_{|\mu|\le 20}
|\Gamma^\mu u'(s,x)|\Bigr)^2 \:dx\:ds 
\end{multline*}
since $\langle c_Ks-r\rangle \gtrsim (1+s)$ on the domain of integration in the left.  Thus, by \eqref{g.1} and \eqref{g.2},
this is also $O(\varepsilon^3)$ for sufficiently small $\varepsilon$, which completes the proof of $(i.)$.

\subsection{High order energy}
Here, we shall show that \eqref{g.2} follows from \eqref{g.1}.
By \eqref{i.4} and \eqref{e.29}, the square of the left side of \eqref{g.2} is 
\begin{multline}
\label{h.1}
\lesssim \varepsilon^2 + \int_0^t \int_\ext \sum_{|\mu|\le 15} |\Gamma^\mu u'| \sum_{|\nu|\le 29} |\Gamma^\nu u''|
\sum_{|\sigma|\le 30} |\Gamma^\sigma u'|\:dx\:ds
\\+\int_0^t \int_\ext \sum_{|\mu|\le 15} |\Gamma^\mu u'| \sum_{|\nu|\le 29} |\Gamma^\nu u''| \frac{\sum_{|\sigma|\le 30}
|\Gamma^\sigma u|}{r}\:dx\:ds
\\+\int_0^t \int_\ext \sum_{|\mu|\le 15} |\Gamma^\mu u''| \Bigl(\sum_{|\nu|\le 30} |\Gamma^\nu u'|\Bigr)^2\:dx\:ds
\\+\int_0^t \int_\ext \sum_{|\mu|\le 15} |\Gamma^\mu u''| \sum_{|\nu|\le 30} |\Gamma^\nu u'| \frac{
\sum_{|\sigma|\le 30} |\Gamma^\sigma u|}{r}\:dx\:ds
\\+\int_0^t \int_\ext \langle x\rangle^{-1} \sum_{|\mu|\le 15} |\Gamma^\mu u'|\Bigl(\sum_{|\nu|\le 30} 
|\Gamma^\nu u'|\Bigr)^2\:dx\:ds
\\+\int_0^t \int_\ext \langle x\rangle^{-1} \sum_{|\mu|\le 15} |\Gamma^\mu u'| \sum_{|\nu|\le 30} |\Gamma^\nu
u'| \frac{\sum_{|\sigma|\le 30} |\Gamma^\sigma u|}{r}\:dx\:ds
\\+\sum_{|\mu|\le 15} \sum_{|\nu|\le 30} \| |\Gamma^\mu u'||\Gamma^\nu u'|\|^2_{L^2_tL^2_x(S_t)}
+\sum_{|\mu|\le 15} \sum_{|\nu|\le 30} \||\Gamma^\mu u'| |\Gamma^\nu u'|(t,\cd)\|_2^2.
\end{multline}
The last eight terms in \eqref{h.1} will be referred to as $I, II,\dots, VIII$ respectively.  Terms
$I$ and $III$ are the key terms; the others are technical terms that result from our
analysis of the perturbed KSS estimates.

We start by bounding the terms $I$ and $II$.  When the spatial integrals are over $|x|\ge c_0 s/2$, 
it follows from \eqref{d.4}, the Schwarz inequality, and a Hardy inequality that these terms are
$$\lesssim \int_0^t \frac{1}{1+s} \sum_{|\mu|\le 17} \|\Gamma^\mu u'(s,\cd)\|_2 
\sum_{|\nu|\le 30} \|\Gamma^\nu u'(s,\cd)\|_2^2\:ds.$$
When the spatial integrals are instead over $|x|< c_0s/2$, we have that $\langle c_Ks-r\rangle\gtrsim
(1+s)$ for any $K=1,\dots, D$.  Thus, by \eqref{d.4}, \eqref{d.8}, and \eqref{g.4}, these terms are
\begin{align*}
&\lesssim \int_0^t \frac{1}{1+s} \sum_{|\mu|\le 17} \|\langle x\rangle^{-1} \Gamma^\mu u'(s,\cd)\|_2
\\&\quad\quad\quad\quad\quad\quad
\times\Bigl(\sum_{|\nu|\le 30} \|\Gamma^\nu u'(s,\cd)\|_2  + s\sum_{|\nu|\le 30} \|\Gamma^\nu u'(t,\cd)
\|_{L^2(|x|<1)} \\&\quad\quad\quad\quad\quad
+ s\sum_{|\nu|\le 16} \|\Gamma^\nu u'(s,\cd)\|_{L^2(|x|<1)}\sum_{|\nu|\le 30} \|\Gamma^\nu u'(s,\cd)\|_2
\Bigr) \sum_{|\sigma|\le 30} \|\Gamma^\sigma u'(s,\cd)\|_2\:ds\\
&\lesssim \varepsilon \int_0^t \frac{1}{1+s} \sum_{|\mu|\le 30} \|\Gamma^\mu u'(s,\cd)\|^2_2\:ds
\\&\quad
+\varepsilon \sum_{|\nu|\le 30} \|\Gamma^\nu u'\|_{L^2_tL^2_x([0,t]\times \{|x|<1\})}\sup_{0\le s
\le t}\sum_{|\sigma|\le 30}
\|\Gamma^\sigma u'(s,\cd)\|_2
\\&\quad\quad\quad\quad\quad 
+ \varepsilon^2 \sup_{0\le s\le t} \sum_{|\nu|\le 30} \|\Gamma^\nu u'(s,\cd)\|^2_2. 
\end{align*}
Here, we have applied the Schwarz inequality and \eqref{g.1}.
Thus, it follows that 
\begin{multline}
\label{h.2}
I+II\lesssim \varepsilon \int_0^t \frac{1}{1+s} \sum_{|\mu|\le 30} \|\Gamma^\mu u'(s,\cd)\|^2_2\:ds
+\varepsilon \sup_{0\le s\le t} \sum_{|\mu|\le 30} \|\Gamma^\mu u'(s,\cd)\|^2_2 \\+ \varepsilon
\sum_{|\mu|\le 30} \|\langle x\rangle^{-5/8} \Gamma^\mu u'\|^2_{L^2_tL^2_x(S_t)}. 
\end{multline}
The last two terms can be bootstrapped and absorbed into the left side of \eqref{h.1}

The bound for the next two terms in \eqref{h.1} is similar.  Again, in $|x|\ge c_0s/2$,
by an application of \eqref{d.4} and \eqref{g.1}, we have that these terms are
$$\lesssim \varepsilon \int_0^t \frac{1}{1+s} \sum_{|\mu|\le 30} \|\Gamma^\mu u'(s,\cd)\|^2_2\:ds.$$
When the inner integrals are over $|x|<c_0s/2$, we may apply \eqref{d.9} and \eqref{g.3} to see that
these terms are
\begin{multline*}
\lesssim \int_0^t \frac{1}{1+s} \sum_{|\mu|\le 18} \|\Gamma^\mu u'(s,\cd)\|_2 \sum_{|\nu|\le 30} \|\Gamma^\nu u'(s,\cd)\|_2^2
\:ds \\+ \sum_{|\mu|\le 15} \|\Gamma^\mu u'\|_{L^2_tL^2_x(S_t\cap\{|x|<1\})}\sum_{|\nu|\le 30}
\|\langle x\rangle^{-1} \Gamma^\nu u'\|_{L^2_tL^2_x(S_t)} \sup_{0\le s\le t} \sum_{|\sigma|\le 30}
\|\Gamma^\sigma u'(s,\cd)\|_2. 
\end{multline*}
Here, we have also applied the Schwarz inequality and a Hardy inequality.  Thus, it follows from \eqref{g.1} 
that
\begin{multline}
\label{h.3}
III+IV \lesssim \varepsilon \int_0^t \frac{1}{1+s} \sum_{|\mu|\le 30} \|\Gamma^\mu u'(s,\cd)\|^2_2\:ds
+\varepsilon \sup_{0\le s\le t} \sum_{|\mu|\le 30} \|\Gamma^\mu u'(s,\cd)\|^2_2 \\+ \varepsilon
\sum_{|\mu|\le 30} \|\langle x\rangle^{-5/8} \Gamma^\mu u'\|^2_{L^2_tL^2_x(S_t)},
\end{multline}
and again, the last two terms will be bootstrapped.

By the Schwarz inequality, a Hardy inequality, \eqref{d.4}, and \eqref{g.1}, we easily obtain
\begin{equation}
\label{h.4}
V+VI\lesssim \varepsilon \sum_{|\mu|\le 30} \|\langle x\rangle^{-1} \Gamma^\mu u'\|_{L^2_tL^2_x(S_t)}
\sup_{0\le s\le t} \sum_{|\nu|\le 30} \|\Gamma^\nu u'(s,\cd)\|_2,
\end{equation} 
which will also be bootstrapped.  Using \eqref{d.4}, the last two terms are also easily handled, resulting in
\begin{equation}
\label{h.5}
VII+VIII\lesssim \varepsilon^2 \Bigl(\sum_{|\mu|\le 30} \|\langle x\rangle^{-1} \Gamma^\mu u'\|^2_{L^2_tL^2_x(S_t)}
+ \sum_{|\mu|\le 30} \|\Gamma^\mu u'(t,\cd)\|^2_2\Bigr).
\end{equation}

If we use the estimates \eqref{h.2}-\eqref{h.5} in \eqref{h.1} and bootstrap the appropriate terms, it follows that
\begin{multline}
\label{h.6}
\sum_{|\mu|\le 30}\|\Gamma^\mu u'(t,\cd)\|_2^2 +\sum_{|\mu|\le 30} \|\langle x\rangle^{-5/8} \Gamma^\mu u'\|^2_{L^2_tL^2_x(S_t)}
\\\lesssim \varepsilon^2 + \varepsilon \int_0^t \frac{1}{1+s} \sum_{|\mu|\le 30} \|\Gamma^\mu u'(s,\cd)\|^2_2\:ds.
\end{multline}
The desired bound, \eqref{g.2}, then follows from Gronwall's inequality, which completes the proof.

\bigskip


\begin{thebibliography}{MA}
\bibitem{AY} R. Agemi and K. Yokoyama: {\em The null condition and global existence of solutions
to systems of wave equations with different speeds}.  Advances in Nonlinear Partial Differential Equations
and Stochastics, (1998), 43--86.
\bibitem{Alinhac2} S. Alinhac: {\em On the Morawetz/KSS inequality for the wave equation on a curved
background}, preprint (2005).
\bibitem{Christodoulou} D. Christodoulou: {\em Global solutions of nonlinear hyperbolic equations for
small initial data}.  Comm. Pure Appl. Math. {\bf 39} (1986), 267--282.
\bibitem{Hid} K. Hidano: {\em An elementary proof of global or almost global existence for quasi-linear
wave equations}.  Tohoku Math. J. {\bf 56} (2004), 271--287.
\bibitem{HidY} K. Hidano and K. Yokoyama: {\em A remark on the almost global existence theorems of Keel,
Smith, and Sogge}.  Funkcial. Ekvac. {\bf 48} (2005), 1--34.
\bibitem{JK} F. John and S. Klainerman: {\em Almost global existence to nonlinear wave equations in
three dimensions}.  Comm. Pure Appl. Math. {\bf 37} (1984), 443-455.
\bibitem{H} L. H\"ormander: {\em Lectures on nonlinear hyperbolic equations}, Springer-Verlag, Berlin, 1997.
\bibitem{Kata} S. Katayama: {\em Global existence for a class of systems of nonlinear wave equations in three
space dimensions}.  Chinese Ann. Math. Ser. B, {\bf 25} (2004), 463--482.
\bibitem{Kata2} S. Katayama: {\em Global existence for systems of wave equations with nonresonant nonlinearities and null
forms}. J. Differential Equations {\bf 209} (2005), 140--171.
\bibitem{Kata3} S. Katayama: {\em A remark on systems of nonlinear wave equations with different propagation speeds}, preprint (2003).
\bibitem{KSS} M. Keel, H.
Smith, and C. D. Sogge: {\em Global existence for a
 quasilinear wave equation
outside of star-shaped domains}. J. Funct. Anal. {\bf{189}}
(2002), 155--226.
\bibitem{KSS2} M. Keel, H. Smith, and C. D. Sogge: {\em Almost
global existence for some semilinear wave equations}. J.
D'Analyse {\bf 87} (2002), 265--279.
\bibitem{KSS3}M. Keel, H.
Smith, and C. D. Sogge: {\em Almost global existence for
quasilinear wave equations in three space dimensions}. J. Amer. Math. Soc. {\bf 17} (2004), 109--153.
\bibitem{K2} S. Klainerman: {\em Uniform decay estimates and the Lorentz invariance of the classical wave
equation}.  Comm. Pure Appl. Math. {\bf 38} (1985), 321--332.
\bibitem{K} S. Klainerman: {\em The null condition and global existence to nonlinear wave equations}. Lect. Appl.
Math. {\bf 23} (1986), 293--326.
\bibitem{KS} S. Klainerman and T. Sideris: {\em On almost global existence for
nonrelativistic wave equations in 3d}. Comm. Pure Appl. Math. {\bf
49} (1996), 307--321.
\bibitem{KY} K. Kubota and K. Yokoyama: {\em Global existences of classical solutions to systems of nonlinear
wave equations with different speeds of propagation}.  Japan. J. Math.  {\bf 27} (2001), 113--202.
\bibitem{LMP} P. D. Lax, C. S. Morawetz, and R. S. Phillips:
{\em Exponential decay of solutions of the  wave equation in the
exterior of a star-shaped obstacle}. Comm. Pure Appl. Math. {\bf
16} (1963), 477--486.
\bibitem{Met} J. Metcalfe: {\em Global existence for semilinear wave equations exterior
to nontrapping obstacles}.  Houston J. Math. {\bf 30} (2004), 259--281.
\bibitem{MNS1} J. Metcalfe, M. Nakamura, and C. D. Sogge: {\em Global existence of solutions
to multiple speed systems of quasilinear wave equations in exterior domains}. Forum.
Math. {\bf 17} (2005), 133--168.
\bibitem{MNS2} J. Metcalfe, M. Nakamura, and C. D. Sogge: {\em Global existence of quasilinear, nonrelativistic
wave equations satisfying the null condition}. Japan. J. Math. {\bf 31} (2005), 391--472.
\bibitem{MS1} J. Metcalfe and C. D. Sogge: {\em Hyperbolic trapped rays and global existence
of quasilinear wave equations}.  Invent. Math. {\bf 159} (2005), 75--117.
\bibitem{MS2} J. Metcalfe and C. D. Sogge: {\em Global existence for Dirichlet-wave equations with quadratic
nonlinearities in high dimensions}.  Math. Ann., to appear.
\bibitem{MS3} J. Metcalfe and C.D. Sogge: {\em Long time existence of quasilinear wave equations
exterior to star-shaped obstacles via energy methods}.  SIAM J. Math. Anal. {\bf 38} (2006), 188--209.
\bibitem{M} C. S. Morawetz: {\em The decay of solutions of the exterior
initial-boundary problem for the wave equation}. Comm. Pure Appl.
Math. {\bf 14} (1961), 561--568.
\bibitem{Si} T. Sideris: {\em Nonresonance and global existence of prestressed
nonlinear elastic waves}.
 Ann. of Math. {\bf 151} (2000), 849--874.
\bibitem{Si2} T. Sideris: {\em The null condition and global existence of
nonlinear elastic waves}.
Invent. Math. {\bf 123} (1996), 323--342.
\bibitem{Si3} T. Sideris and S.Y. Tu: {\em Global existence
for systems of nonlinear wave equations in 3D with multiple
speeds}. SIAM J. Math. Anal. {\bf 33} (2001), 477--488.
\bibitem{SS} H. Smith and C. D. Sogge: {\em Global Strichartz estimates for
nontrapping perturbations of the Laplacian}. Comm. Partial
Differential Equations {\bf 25} (2000), 2171--2183.
\bibitem{S}  C. D. Sogge: {\em Lectures on nonlinear wave equations}.
International Press, Cambridge, MA, 1995.
\bibitem{So2} C. D. Sogge: {\em Global existence for nonlinear
wave equations with multiple speeds}.  Harmonic analysis at Mount Holyoke (South Hadley,
MA, 2001), 353--366, Contemp. Math., 320, Amer. Math. Soc., Providence, RI, 2003.
\bibitem{Sterb} J. Sterbenz: {\em Angular regularity and Strichartz estimates
for the wave equation} with an appendix by I. Rodnianski.
Int. Math. Res. Not. {\bf 2005}, 187--231.
\bibitem{Strauss} W. A. Strauss: {\em Dispersal of waves vanishing on the boundary of an exterior domain}.
Comm. Pure Appl. Math. {\bf 28} (1975), 265--278.
\bibitem{Y} K. Yokoyama: {\em Global existence of classical solutions to systems of wave equations
with critical nonlinearity in three space dimensions}.  J. Math. Soc. Japan {\bf 52} (2000), 609--632.
\end{thebibliography}
\end{document}